\documentclass[sn-mathphys,Numbered]{sn-jnl}
\usepackage{graphicx}
\usepackage{multirow}%
\usepackage{amsmath,amssymb,amsfonts}%
\usepackage{amsthm}%
\usepackage{mathrsfs}%
\usepackage[title]{appendix}%
\usepackage{xcolor}%
\usepackage{textcomp}%
\usepackage{manyfoot}%
\usepackage{booktabs}%
\usepackage{algorithm}%
\usepackage{algorithmicx}%
\usepackage{algpseudocode}%
\usepackage{listings}%
\usepackage{comment}
\usepackage{pgfplots} 				
\usepackage{tikz}
\usepackage{tikzscale}
\usepackage{environ}
\usepackage{cleveref} 

\usepackage{tikz}
\usetikzlibrary{shapes,backgrounds,calc}

\usetikzlibrary{decorations.pathreplacing}
\usetikzlibrary{arrows.meta}
\usepackage[caption=false]{subfig}
\definecolor{backcolor}{rgb}{.7,.7,1}
\definecolor{backcolor2}{rgb}{1,.7,0.7}
\usetikzlibrary{positioning}
\makeatletter
\tikzset{circle split part fill/.style  args={#1,#2}{%
 alias=tmp@name, 
  postaction={%
    insert path={
     \pgfextra{%
     \pgfpointdiff{\pgfpointanchor{\pgf@node@name}{center}}%
                  {\pgfpointanchor{\pgf@node@name}{east}}%
     \pgfmathsetmacro\insiderad{\pgf@x}
      \fill[#1] (\pgf@node@name.base) ([xshift=-\pgflinewidth]\pgf@node@name.east) arc
                          (0:180:\insiderad-\pgflinewidth)--cycle;
      \fill[#2] (\pgf@node@name.base) ([xshift=\pgflinewidth]\pgf@node@name.west)  arc
                           (180:360:\insiderad-\pgflinewidth)--cycle;            
         }}}}}  
 \makeatother

\tikzset{
    position/.style args={#1:#2 from #3}{
        at=(#3), anchor=#1+180, shift=(#1:#2)
    }
}

\def\ls{\textcolor{orange}}

\newtheorem{example}{Example}%
\newtheorem{remark}{Remark}%

\newtheorem{definition}{Definition}%
\newcommand{\R}{\mathbb{R}\,}
\begin{document}

\title[Comparison study for high dimensional approximation]{A comparison study of supervised learning techniques for the approximation of high dimensional functions and feedback control}

\author[1]{\fnm{Mathias} \sur{Oster}}\email{oster@igpm.rwth-aachen.de}

\author*[2]{\fnm{Luca} \sur{Saluzzi}}\email{luca.saluzzi@sns.it}

\author[3]{\fnm{Tizian} \sur{Wenzel}}\email{tizian.wenzel@uni-hamburg.de}

\affil[1]{\orgdiv{IGPM}, \orgname{RWTH Aachen}, \orgaddress{\street{Templergraben 55}, \city{Aachen}, \postcode{52062}, \country{Germany}}}

\affil*[2]{\orgdiv{Department of Mathematics}, \orgname{Scuola Normale Superiore}, \orgaddress{\street{P.za dei Cavalieri,
7}, \city{Pisa}, \postcode{56126},  \country{Italy}}}

\affil[3]{\orgdiv{Department of Mathematics}, \orgname{Universität Hamburg}, \orgaddress{\street{Bundesstraße 55}, \city{Hamburg}, \postcode{20146}, \country{Germany}}}

\abstract{Approximation of high dimensional functions is in the focus of machine learning and data-based scientific computing.
In many applications, empirical risk minimisation techniques over nonlinear model classes are employed. 
Neural networks, kernel methods and tensor decomposition techniques are among the most popular model classes.
We provide a numerical study comparing the performance of these methods on various high-dimensional functions with focus on optimal control problems, 
where the collection of the dataset is based on the application of the State-Dependent Riccati Equation.}


\maketitle

\section{Introduction}

Finding surrogates for functions in high dimensions has become one of the key tasks in scientific computing. 
One interesting example is the approximation of the value function of the optimal control of PDEs, which leads, after semi-discretisation, to very high-dimensional control problems. 
In optimal control the value function plays a crucial role as it can provide optimal feedback laws. There are two major ways to calculate the value function. 
First, one could solve the Bellman or Hamilton-Jacobi-Bellman (HJB) equation \cite{bellman1966dynamic, DeterministicHJB, falcone2013semi}.
For both equations there is a wide range of numerical tools in optimal feedback control, see e.g. \cite{pol_approx_kunisch, alla2015efficient,Zhao,linearlikePI,OnlinePolicy,DataHJB,SOS,safetySOS,SemiLagranigian,SemiLagrangianStochastic,Falcone1987,FALCONESplitting,VIM,ALLA2020192,maxplus_det,maxplus_stoch}.
Second, simultaneously to the work of Bellman, Pontryagin developed a set of necessary conditions for optimal problems \cite{Pontryagin}, the so-called Pontryagin maximum principle (PMP).
This enables a pointwise evaluation of the value function and can thus be used for the recovery of the value function \cite{beeler2000feedback, kang2017mitigating, nakamura2019adaptive, azmi2020optimal}.

Most classical function approximation tools suffer from the curse of dimensionality, 
i.e.\ the exponential growth of complexity with respect to the input dimension. 
To break the curse various methods have been developed, most prominently neural networks, kernel methods and tensor decomposition techniques. 
In most applications, these methods are based on empirical Least-Squares methods  on a nonlinear model class, known as empirical risk minimization techniques \cite{vapnik1992principles, steinwart2008support} in statistic learning. 
For example, to approximate a function $f:\Omega\mapsto\mathbb R$ defined on some domain $\Omega \subset \mathbb{R}^d$, 
one minimizes the functional $$\frac 1 N \sum_{i=1}^N |f(x_i)-f_M(x_i)|^2$$
for $N$ samples $x_i \in \Omega$ distributed according to some density $\rho$ within some model class $f_M\in\mathcal M\subset C(\Omega)$. 
This means that one needs to be able to access the target function $f$ on given or chosen samples. 
Due to the high dimensionality of the input space $\Omega$, one usually resorts to nonlinear model classes providing improved expressibility compared to linear models of the same complexity, at the cost of increasingly difficult optimization tasks when minimizing the empirical risk.

One of these nonlinear model classes is using structured representations of polynomials like hierarchical tensor formats,
which allow to reduce the number of parameters within the coefficient tensor of a linear ansatz \cite{hackbusch-2012}.
More precisely, we consider a sub-manifold in $\otimes_{j=1}^d \mathbb{R}^{n_i}$
defined by multi-linear parametrizations.  
Here we use tensor trains which are a special case of a hierarchical or tree based tensor format \cite{hackbusch-2012}. 
Tensor trains have been invented by \cite{Oseledets, osel-tt-2011} and applied to various high-dimensional PDE's \cite{Khoromskij-book}, 
however the parametrization has already been already used in quantum physics much earlier.
For good  surveys  we refer to 
\cite{Hackbusch2014,Bachmayr-Uschmajew-Schneider,Legeza-Schneider,Hackbusch-Acta}. 
The tensor train representation has appealing properties, making them attractive for treatment of the present problems. 
For example they contain sparse polynomials, 
but are much more flexible at a price of a slightly larger overhead, see e.g.\ \cite{Dahmen-Bachmayr} for a comparison concerning parametric PDEs. 

There has been extensive use of tensor trains in high dimensional optimal control problems \cite{DKK21,oster2,oster3}
and stochastic control problems \cite{stefansson,horowitz2014linear, oster1, gorodetsky2018high}. 
In particular, in this paper we will focus on two specific tensor train approaches: 
the TT Gradient Cross \cite{dolgov2022data} and the block-sparse tensor train \cite{Micha}. 
Cross approximation methods \cite{ot-ttcross-2010,so-dmrgi-2011proc,grasedyck-par-cross-2015,sav-qott-2014} have been introduced in order to adjust the sampling sets to reduce the conditioning of the interpolation problem and enhance the accuracy of approximation. 
TT Gradient Cross makes use of the Cross interpolation to construct efficiently the interpolation indices and takes into account the information of the gradient of the target function to improve the stability. 
The block-sparse tensor train approach is exploiting sparsity patterns in the cores of tensor trains to avoid overparametrization by identifying homogeneous degree basis functions.

Another approach is the use of kernel methods \cite{wendland2005scattered},
which comprise various techniques in numerical approximation, machine learning and scientific computing.
These flexible tools allow to work with arbitrarily scattered data in high dimensional space and allow for a convenient mathematical analysis based on reproducing kernel Hilbert spaces.

Probably the most popular and widespread techniques in supervised learning are neural networks. 
The last decade has seen an incredible development in machine learning since efficient back-propagation and optimization algorithms as well as powerful hardware allowed for very complicated neural network architectures, 
see \cite{DeepLearningGitta,DHP,DeepLearningWeinan,DeepLearningHigham} for an introduction to neural networks from a mathematical perspective. 
They are employed in virtually all machine learning and scientific computing tasks \cite{ImageRec,DeepLearningPDE}. 
Also in optimal control they were used to obtain (sub)optimal feedback laws and surrogates for the value functions \cite{Kunisch,Kunisch2,darbon2020overcoming, nusken2020solving,ito2020neural,DEMO2023383,Han_Jentzen_E_2018,sympocnet,Zhou_2021,Onken2021,ruthotto2020machine,ABK21,grune2020computing}.

Lastly, we like to mention that also sparse polynomial techniques have been successfully employed in the context of high dimensional optimal control \cite{Kunisch3, Kunisch4, azmi2021optimal}. 
However, they are not part of the numerical comparison in this study. \\



In this paper we conduct a numerical study comparing Tensor Trains, neural networks and kernel methods by approximating high dimensional functions and applying them to value functions from high dimensional optimal control problems. 
After introducing the general optimal control problem in Section \ref{sec:OC}, we recall our nonlinear model classes in Section \ref{sec:ML}. 
In Section \ref{sec:NE} we provide a variety of numerical examples where we compare the performance of the different methods. 
Section \ref{sec:conclusion} concludes the paper.

\section{The Optimal Control Problem}\label{sec:OC}

In this section we present the infinite horizon problem and the corresponding synthesis of the feedback law. 
We first introduce the optimal control framework using the Hamilton-Jacobi-Bellman (HJB) formalism, 
then we pass to the synthesis of suboptimal feedback laws which will provide the data for the supervised learning techniques.
We consider a dynamical system in control affine-form given by
\begin{equation}\label{eq}
\left\{ \begin{array}{l}
\dot{y}(s)=f(y(s)) + B(y(s))u(s), \;\; s\in(0,+\infty),\\
y(0)=x\in\mathbb{R}^d.
\end{array} \right.
\end{equation}
We denote by $y:[0,+\infty)\rightarrow\R^d$ the state of the system, by $u:[0,+\infty)\rightarrow\R^m$ the control signal and by $\mathcal{U}=L^\infty ([0,+\infty);U)$ the set of admissible controls where $U\subset \R^m$. The system dynamics $f:\R^d\rightarrow\R^d$ and $B:\R^d\rightarrow\R^d$ are assumed to be $\mathcal{C}^1(\R^d)$ functions.

\noindent We introduce the infinite horizon cost functional:
\begin{equation}\label{cost}
 J(u(\cdot,x)):=\int_0^{+\infty} r(y(s))+ u^{\top}(s) R u(s)\,ds\,,
\end{equation}
where $r:\R^d\rightarrow\R^+$ and $R\in \R^{m\times m}$ is a symmetric positive definite matrix. Our aim is to compute an optimal control in feedback form, $e.g.$  a control signal fully determined upon the current state of the system. We start by the definition of the value function for a given initial condition $x \in\R^d$:
\begin{equation}
V(x):=\inf\limits_{u\in\mathcal{U}} J(u(\cdot,x))\,,
\label{VF}
\end{equation}
which satisfies the following Hamilton-Jacobi-Bellman PDE for every $x\in\R^d$
\begin{equation}\label{HJB}
\min\limits_{u\in U }\left\{(f(x)+B(x)u)^{\top} \nabla V(x) + r(x)+ u^{\top} R u \right\}=0.
\end{equation}
The HJB PDE \eqref{HJB} is a first-order fully nonlinear PDE defined over $\R^d$, where $d$ represents the dimension of the considered dynamical system. The dimension of the problem may be large and this limitation is known as curse of dimensionality. In the case, $i.e.$ $U=\R^m$, the minimization in the equation \eqref{HJB} can be computed explicitly as
\begin{equation}\label{optc}
u^*(x)=-\frac{1}{2} R^{-1}B(x)^{\top} \nabla V(x)\,,
\end{equation}
leading to the following version of the HJB PDE
	\begin{align}\label{hjbu}
	\nabla V(x)^{\top}f(x)-\frac14\nabla V(x)^{\top}B(x)R^{-1}B(x)^{\top}\!\nabla V(x)\!+r(x)=0\,.
\end{align}

In this work, instead of considering directly the high-dimensional HJB PDE \eqref{hjbu}, 
we are interested in retrieving an approximation of the value function using supervised learning techniques. 
Indeed, we will consider an approximation of $V(x)$ in a regression framework, where we assume measurements of the function at sampling points. 

We are going to consider a fast, but suboptimal alternative: the State-Dependent Riccati Equation (SDRE). 
This will allow to generate a synthetic dataset that allows to approximate the value function, leading to the synthesis of a suboptimal control which is still able to stabilize the dynamical system.

\subsection{State-Dependent Riccati Equation}

The State Dependent Riccati Equation (SDRE) is a powerful mathematical tool that finds widespread applications \cite{ccimen2008state,allasdre}. Originating from the classical Riccati Equation, the SDRE extends its utility by incorporating state-dependent coefficients, thereby accommodating systems with nonlinearity and time-varying dynamics. The approach relies on sequential resolution of linear-quadratic control problems that stem from progressive linearization of the dynamics along a trajectory.

Let us suppose that the cost functional \eqref{cost} can be rewritten in the following form
\begin{equation}
    J_{\infty}(u;x) = \int\limits_0^{+\infty} y(t)^\top Q(y) y(t) + u(t)^\top R u(t)\, dt \,,
    \label{quadratic_cost}
\end{equation}
with $Q: \R^d \rightarrow \mathbb{R}^{d\times d}, Q(y)\succeq 0$ $\forall y \in \R^d$ and the dynamical system expressed in semilinear form
\begin{align}
    \dot{y}(t) & = A(y(t)) y(t) +B(y(t)) u(t) \\
    y(0) & = x\,.
    \label{semilinear}
\end{align}
Assume the dynamics is linear in the state, $i.e.$ $A(y(t)) =A\in\mathbb{R}^{d\times d}$ and $B(y(t))=B\in\mathbb{R}^{d\times m}$, and that the matrix $Q(y)=Q$ is constant in the state. This is called Linear Quadratic Regulator (LQR) problem. If the pair $(A,B)$ is stabilizable and the pair $(A,Q^{1/2})$ is detectable, then the optimal feedback control is computed via the following formula
\begin{equation}
   u(y) = -R^{-1} B^{\top} Py.
\label{control_SDRE}
\end{equation}
Here, $P\in\mathbb{R}^{d\times d}$ is the unique positive definite solution of the Algebraic Riccati Equation (ARE)
\begin{align*}
A^\top P + P A
  -PBR^{-1}B^\top P +Q = 0\,.
\end{align*}
 Essentially, the SDRE technique extends this approach introducing the dependence on the current state, that is,
\begin{equation}
   u(y) = -R^{-1} B^{\top}(y) P(y)y\,,
    \label{control_sdre_feed}
\end{equation}
where $P(y)$ now is the solution of a state-dependent ARE
\begin{align}
A^\top(y) P(y) + P(y) A(y)-P(y)B(y)R^{-1}B^\top(y)P(y)+Q & = -Q(y),
    \label{sdre}
\end{align}
where $A(y)$ and $B(y)$ are fixed at the state $y$.
This procedure can be iterated along the trajectory, solving sequentially \eqref{sdre} as the state $y(t)$ is evolving in time.

Assuming suitable stability hypothesis, it is possible to show that the closed loop dynamics generated by the feedback law \eqref{control_sdre_feed} are locally
asymptotically stable (we refer to \cite{Banks_Lewis_Tran_2007} for the exact statement and further details).

\section{Machine learning methods}\label{sec:ML}

\subsection{Machine learning and least-squares loss}
In this article we focus on the use of supervised machine learning techniques for high-dimensional functions and optimal control problems. 
In each of the following numerical example we aim to approximate a target function $f$. 
One way to do so is to minimize an empirical $L^2$ loss between the target function and a nonlinear model class, 
i.e. we are considering a problem of the form 
\begin{equation}\label{eq:mse_loss}
    \min_{f_M\in\mathcal M}\frac 1 N \sum_{j=1}^N |f(x_j)-f_M(x_j)|^2,
\end{equation}
where $x_j$ are samples in $\Omega$ and $\mathcal M$ is the model class.
This allows us to use random samples. 
Another route is taken by the TT-Cross algorithm introduced in the following. 
There an active learning strategy is employed to reduce the sample complexity.

All methods which will be introduced in the following do employ a global function approximation, which is in contrast to localized methods such as as finite elements. 
Furthermore, all model classes will be nonlinear except the kernel method. 



\subsection{Tree Based Tensor Representation - Tensor Trains}\label{sect:TT}

For the approximation of the value function \eqref{VF}, we define a nonlinear model class to circumvent the curse of dimensionality. 
To this end, we choose an underlying finite dimensional subspace for the approximation of the sought value function. For the present purpose we 
take a space $\Pi_{i,n_i} = \mathrm{span} \{\psi_{i_1}, \dots, \psi_{i_d} \}$ of one-dimensional polynomials of degree smaller than $n_i$  and consider the tensor product of such polynomial spaces
$$ \mathcal{V}_p := \Pi_{1,n_1} \otimes \cdots \otimes \Pi_{d,n_d}.
$$
This is a space of multivariate (tensor product) polynomials with bounded multi-degree.
Its elements $v \in \mathcal V_p$ can be represented as
\begin{equation}\label{eq:function_fulltensor}
v ( x_1 , \ldots , x_d ) = 
\sum_{i_1 , \ldots , i_d = 0}^{n_1, \dots, n_d}
c_{ i_1, \ldots , i_d} \psi_{i_1} (x_1) \cdots 
\psi_{i_d} (x_d),
\end{equation}
exhibiting that $c\in \mathbb R^{n_1, \dots, n_d}$ suffers from the curse of dimensionality.
For the sake of readability we will henceforth write $c[i_1, \dots, i_d] = c_{i_1, \dots, i_d}$ and say that $c$ is an order $d$ tensor.

The tensor train decomposition aims to represent an order $d$ tensor by a sequence of order $3$ tensors, connected by contractions.
This means that we represent $c$ by $U_1 \in \mathbb R^{n_1, r_1}$, $U_2 \in \mathbb R^{r_1, n_2, r_2}, \dots, U_{d-1} \in \mathbb R^{r_{d-2}, n_{d-1}, r_{d-1}}$ and $U_d \in \mathbb R^{r_{d-1}, n_d}$ such that
\begin{equation}\label{eq:TT_representation_coefficienttensor}
    c[i_1, \dots, i_d] = \sum_{j_1 = 1}^{r_1} \dots \sum_{j_{d-1} = 1}^{r_{d-1}} U_1[i_1, j_1] U_2 [j_1, i_2, j_2] \dots U_d[j_{d-1}, i_d].
\end{equation}
The TT-rank is introduced as the element wise smallest tuple $\mathbf r = (r_1, \dots, r_{d-1})$ such that a decomposition of the form \eqref{eq:TT_representation_coefficienttensor} exists.
The TT-rank is well defined and,
denoting $r = \max \{r_i\}$ and $n = \max \{n_i\}$, the tensors of fixed TT-rank form a smooth manifold of dimension $\mathcal O (dnr^2)$, which means that for fixed ranks the dimension of the manifold does increase linearly with the order $d$. Quadratic functions \cite{DKK21} or weakly correlated Gaussian functions \cite{rdgs-tt-gauss-2022}, for example, admit an approximation with TT-rank growing at most polynomial in $d$ and poly-logarithmic in the approximation error.

\begin{figure}[h!]
    \centering
    \begin{tikzpicture}
\begin{scope}[every node/.style={scale=1,minimum size=5mm, draw,  fill=white}]
    \node (A1) at (0,0) {$U_1$}; 
    \node (A2) at (2,0) {$U_2$}; 
    \node (A3) at (4,0) {$U_3$}; 
    \node (A4) at (6,0) {$U_4$}; 
	\node [ position=-90:1 from A1](B1) {$\Psi_1(x_1)$};
   	\node [ position=-90:1 from A2](B2) {$\Psi_2(x_2)$};
   	\node [ position=-90:1 from A3](B3) {$\Psi_3(x_3)$};
   	\node [ position=-90:1 from A4](B4) {$\Psi_4(x_4)$};

\end{scope}
    \node (C) at (-4,0) {$v(x)$}; 
    \node (D) at (-2,0) {$=$}; 

\begin{scope}[every edge/.style={draw=black,thick}]
	\path [-] (A1) edge node[midway,left] [above] {$r_1$} (A2);
	\path [-] (A3) edge node[midway,left] [above] {$r_2$} (A2);
	\path [-] (A3) edge node[midway,left] [above] {$r_3$} (A4);
	
	\path [-] (A1) edge node[midway,left]  {$n_1$} (B1);
	\path [-] (A2) edge node[midway,left]  {$n_2$} (B2);  
	\path [-] (A3) edge node[midway,left]  {$n_3$} (B3);  
	\path [-] (A4) edge node[midway,left]  {$n_4$} (B4);  
\end{scope} 
              
\end{tikzpicture}
    \caption{Graphical representation of TT tensor train induced polynomial in four variables.}
    \label{TT_polynomial}
\end{figure}
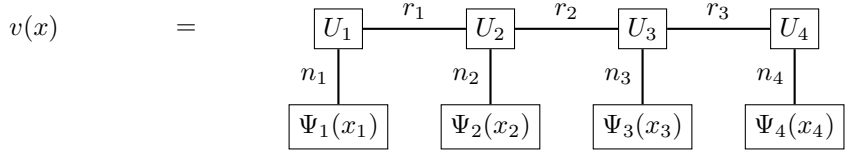

Observing that the component tensors $U_i$ are connected via a single contraction/summation to $U_{i-1}$ and $U_{i+1}$, we can represent the decomposition in a graph, by setting the components $U_i$ as nodes and indicate contractions by links between the nodes.
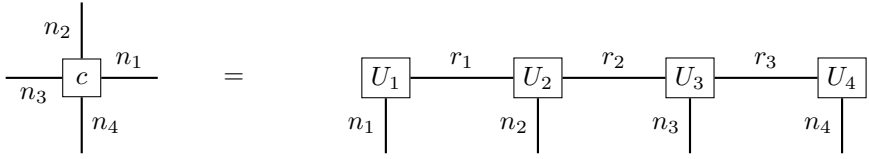
\begin{figure}[h!]
    \centering
    \begin{tikzpicture}
\begin{scope}[every node/.style={scale=1,minimum size=5mm, draw,  fill=white}]
    \node (A1) at (0,0) {$U_1$}; 
    \node (A2) at (2,0) {$U_2$}; 
    \node (A3) at (4,0) {$U_3$}; 
    \node (A4) at (6,0) {$U_4$}; 
    \node (C) at (-4,0) {$c$}; 
\end{scope}
    \node (D) at (-2,0) {$=$}; 
\begin{scope}[every edge/.style={draw=black,thick}]
	\path [-] (A1) edge node[midway,left] [above] {$r_1$} (A2);
	\path [-] (A3) edge node[midway,left] [above] {$r_2$} (A2);
	\path [-] (A3) edge node[midway,left] [above] {$r_3$} (A4);
	\path [-] (A1) edge node[midway,left] {$n_1$} (0,-1);
	\path [-] (A2) edge node[midway,left] {$n_2$} (2,-1);
	\path [-] (A3) edge node[midway,left] {$n_3$} (4,-1);
	\path [-] (A4) edge node[midway,left] {$n_4$} (6,-1);
	\path [-] (C) edge node[midway,left] [above] {$n_1$} (-3,0);
	\path [-] (C) edge node[midway,left] {$n_2$} (-4,1);
	\path [-] (C) edge node[midway,left] [below] {$n_3$} (-5,0);
	\path [-] (C) edge node[midway,left] [right]{$n_4$} (-4,-1);
\end{scope} 
\end{tikzpicture}
    \caption{Graphical representation of a TT representation of $c$ in four variables.}
\end{figure}
In the next step we plug the TT-decomposition of the coefficient tensor \eqref{eq:TT_representation_coefficienttensor} into the representation in \eqref{eq:function_fulltensor}.
    To this end we introduce the short form $\Psi_i(x_i) = [\psi_{i_1}(x_i), \dots, \psi_{i_d}(x_i)] \in \mathbb R^{n_i}$.
Then 
\begin{multline}
    v(x_1,\dots,x_d) =\sum_{i_1,\dots,i_d}^{n_1,\dots,n_d} \sum_{j_1,\dots,j_{n-1}}^{r_1,\dots,r_{n-1}} U_1[i_1,j_1] U_2[j_1, i_2, j_2] \dots U_d[j_{d-1},i_d] \\
    \big(\Psi_1(x_1)\big)[i_1](\Psi_2(x_2)\big)[i_2]\cdots \big(\Psi_d(x_d)\big)[i_d],
    \label{FTT}
\end{multline}
which means that every open index of the TT representation is contracted with the one-dimensional basis functions. This representation is known as Functional Tensor Train (FTT) format of the function $v$.
The graphical representation of this tensor network is given in Figure \ref{TT_polynomial}.
Note that any basis can be chosen for $\Psi_i$. 
In this paper we use a set of orthonormal polynomials.
In this case,  we have a Parseval formula providing a norm equivalence between the function space and the Frobenius norm of the coefficients, which guarantees stability of our representations. 

It turns out, that optimization procedures in this TT format can be solved by consecutively optimizing one component $U_l$ while the others are fixed. This alternating Least-Squares  (ALS) algorithm converges at least to a local minimum \cite{ALS}.

\subsubsection{Block-Sparse Tensor Trains}

One important observation is that in general, tensor trains will parameterize polynomials with high-mixed degree which might lead to numerical instabilities if the sought function has a bounded maximal degree. To overcome this superfluous degrees of freedom one can employ the so-called block sparse tensor trains \cite{Micha}. 

As it turns out, homogeneous polynomials of degree $\tilde g$ exhibit a representation as tensor trains for which the core admit a sparse representation with block sizes $\rho_{k,\tilde g}$ which also provide rank bounds. 
Very importantly these block sparse structures are preserved under essential tensor train manipulation as TT-SVD and rounding. 
Furthermore, the ALS algorithm can be restricted to respect these sparsity patterns. 
By introducing an extra index in the last core one can also parameterize non-homogeneous polynomials in a block sparse fashion.

 Let us give a small example on how the cores will look like for the block-sparse TTs.
\begin{example}[\textit{Block Sparsity}]
    Let $p=4$ and $g=3$ be given and let $c$ be a tensor train such that $Lc = gc$.
    Then for $k=2,\ldots,d-1$ the component tensors $C_k$ of $c$ exhibit the following block sparsity (up to permutation).
    For indices $i$ of order $r_{k-1}$ and $j$ of order $r_k$
    \begin{align*}
    C_k(i,1,j) =& \begin{pmatrix}* &0&0&0\\0&*&0&0\\0&0&*&0\\0&0&0&* \end{pmatrix}\quad C_k(i,2,j) = \begin{pmatrix}0 &*&0&0\\0&0&*&0\\0&0&0&*\\0&0&0&0 \end{pmatrix}\\
    C_k(i,3,j) =& \begin{pmatrix}0 &0&*&0\\0&0&0&*\\0&0&0&0\\0&0&0&0 \end{pmatrix}\quad
    C_k(i,4,j) = \begin{pmatrix}0 &0&0&*\\0&0&0&0\\0&0&0&0\\0&0&0&0 \end{pmatrix}.
    \end{align*}
\end{example}

Another structural assumption that frequently can be found in application is a certain grouping of variables, 
i.e. the function can be written as a sum of sub-functions each of which depend only on few variables. 
To make use of this, the following notion of locality is introduced:
\begin{definition}[\cite{Micha}]
    Let $u\in W_g^d$ be a homogeneous polynomial and $B$ be the symmetric coefficient tensor.
    We say that $u$ has a variable locality of $K_{\mathrm{loc}}$ if $B(\ell_1,\ldots,\ell_g)=0$ for all $(\ell_1,\ldots,\ell_g)\in\mathbb{N}_d^g$ with
    \[
        \max\{|\ell_{m_1}-\ell_{m_2}|\,:\,m_1,m_2=1,\ldots,g\} > K_{\mathrm{loc}} .
    \]
\end{definition}

Note that the locality and degree bound will give rank bounds.
All in all, we will employ this block-sparse tensor train to reduce the complexity of the ansatz space. 
We will comment in the numerical examples, where this ansatz will be too restrictive and where it can provide benefits to the usual TTs.

\subsubsection{TT-Gradient Cross}

Given a target function $V$ and its FTT representation \eqref{FTT} $\tilde{V}$, the TT-Gradient Cross algorithm aims to solve for a given sample points $\{x_i \}_{i=1}^N$ and a dataset $\{ \,V(x_i),\, \nabla V(x_i)\}_{i=1}^N$ the following regression problem

\begin{equation}
\min_{U_1, \ldots, U_d} \sum_{i=1}^N | \tilde{V}(x_i) - V(x_i)|^2 + \lambda \Vert \nabla \tilde{V}(x_i) -\nabla V(x_i) \Vert^2,
\label{regr_gradient}
\end{equation}
where $\lambda$ is a parameter tuning the gradient information. 
The problem can be attacked using fast algorithms based on the so-called \emph{cross interpolation} \cite{ot-ttcross-2010}. 
Given a prefixed set of collocation points $X_1 \times \ldots \times X_d$, we apply an alternating direction strategy solving sequentially least square problems. 
At the $k-$th iteration our goal is to find interpolation sets $\overline{X}_{<k} \subset X_1 \times \cdots \times X_{k-1}$ and $\overline{X}_{>k} \subset X_{k+1} \times \cdots \times X_{d}$ with $r_{k-1}$ and $r_k$ points, respectively.
Let us suppose that in the $k$-th step the sets $\overline{X}_{<k}$ and $\overline{X}_{>k}$ are given.
We point out that the number of unknowns in $U_k$ and the cardinality of the set $\overline{X}_{<k} \oplus X_k \oplus \overline{X}_{>k} $ is equal to $r_{k-1} n_k r_k$.
Solving the least square problem related to actual sampling points, 
one can compute the current $U_k$ and thanks to pivoting techniques (for example the $maxvol$ method \cite{gostz-maxvol-2010}), 
it is possible to select the next sampling sets $X_{<k+1}$ and $X_{>k-1}$ as subsets of $X_{<k} \oplus X_k$ and $X_k \oplus X_{>k}$.
This step can be iterated for all $k=1,\ldots,d$ and the TT cores are updated accordingly until the algorithm converges. 
In this case we claim that the method converges if the norm of the difference of the coefficients in Frobenius norm computed in two consecutive steps is below a certain threshold denoted as $tol_{stop}$.

In the context of optimal control problems, the target function $V$ is the value function and the surrogate model $\tilde{V}$ helps  for a fast synthesis of feedback controls. Indeed, the optimal control is given by
$$
u(x)=-\frac{1}{2}R^{-1} B(x)^{\top} \nabla V(x),
$$
 where $R$ and $B(x)$ have been introduced in Section \ref{sec:OC}.
 Computed an approximation of the FTT representation \eqref{FTT} by the TT Gradient Cross, in the online phase we need to compute the gradient of the surrogate model obtaining the feedback control
 $$
\tilde{u}(x)=-\frac{1}{2}R^{-1} B(x)^{\top} \nabla \tilde{V}(x).
$$
The TT format enables to compute the gradient in $O(dnr^2)$ operations per point, resulting in an expedited synthesis of the feedback control.
For a detailed description of the method and its application we refer to \cite{dolgov2022data}.

\subsection{Kernel methods}

Another popular method in machine learning are summarized as kernel methods. 
This class of methods revolves around the use of a kernel $k$, 
which is a symmetric function $k: \Omega \times \Omega \rightarrow \mathbb{R}$, 
that satisfies some definiteness properties like strict positive definiteness, 
i.e.\ the kernel matrix $(k(x_i, x_j))_{i,j=1}^n$ is positive definite for any choice of pairwise distinct points $\{ x_i \}_{i=1}^n$ and any $n \in \mathbb{N}$.
Some popular examples are the Gaussian kernel or the exponential kernel,
    \begin{align}
    \label{eq:kernels_example}
        \begin{aligned}
        k_\text{Gaussian}(x, y) &= \exp(-\Vert x - y \Vert^2), \\
        k_\text{exp}(x, y) &= \exp(-\Vert x - y \Vert)
        \end{aligned}
    \end{align}
which are radial basis function kernels in $\mathbb{R}^d$ for any $d \in \mathbb{N}$.
It turns out, that every strictly positive definite kernels gives rise to a unique reproducing kernel Hilbert spaces (RKHS) $\mathcal{H}_k(\Omega)$, 
which allows for a thorough theoretical analysis.
Under mild assumptions on the domain $\Omega$, e.g.\ a Lipschitz boundary, these RKHS can frequently characterized in terms of Sobolev spaces.
For example, the RKHS of the exponential kernel from Eq.~\eqref{eq:kernels_example} is norm equivalent to the Sobolev space $H^{(d+1)/2}(\Omega)$,
while the RKHS of the Gaussian kernel consists of analytic functions and a full characterization is more sophisticated. 
A representer theorem for kernel approximation \cite{steinwart2008support} states that the optimal solution for the MSE loss task of Eq.~\eqref{eq:mse_loss} can be found as
    \begin{align}
    \label{eq:kernel_model}
        s_X(\cdot) = \sum_{j=1}^M \alpha_j k(\cdot, x_j),
    \end{align}
whereby the coefficients $\{ \alpha_j \}_{j=1}^M$ can be frequently computed directly.
Such kernel models are used for statistical learning \cite{steinwart2008support}, 
numerical approximation \cite{wendland2005scattered}, 
PDE approximation \cite{chen2021solving} 
and machine learning \cite{meanti2022efficient}, 
among others.
Recent machine learning research also aims at modifying the kernel, in order to obtain data-adapted or deep kernel models \cite{owhadi2019kernel, suykens2017deep, wenzel2023data}. 
In numerical approximation, 
based on assumptions like $f \in \mathcal{H}_k(\Omega)$ or also weaker ones,
sharp error estimates for the residual $f - s_X$ on some Lipschitz domain $\Omega \subset \mathbb{R}^d$ can be derived in various norms, e.g.\ $\Vert \cdot \Vert_{L^2(\Omega)}, \Vert \cdot \Vert_{L^\infty(\Omega)}$ or even Sobolev norms \cite{narcowich2005sobolev, wendland2005approximate}. 
Similar to other approximation methods, such error estimates suffer the curse of dimensionality,
however target data adapted point choices $X \subset \Omega$ offer the possibility to break it \cite{wenzel2023analysis}. 
Though the computation of the optimal weights $\{ \alpha_j \}_{j=1}^M$ can be done explicitly, this can be challenging from the computational point of view for large amount of data points $M \gg 1$.
Thus recent research aims at scaling kernel methods to large amounts of data, e.g.\ via iterative preconditioned training methods and efficient implementations \cite{meanti2022efficient, ma2019kernel}. 


\subsection{Neural Networks}

The most prominent machine learning tools are (deep) neural networks \cite{goodfellow2016deep},
that are used especially for high-dimensional approximation tasks like image recognition or image generation. 
Although NN techniques proved to be very powerful for real world applications, it turns out that they pose difficulties to obtain quantitative results for approximation theoretical bounds due to the iterative and non-convex training procedures.
There is a plethora of sophisticated architectures and structures available like convolutional neural networks, residual neural network and recurrent neural networks aiming on different aspects of machine learning tasks as image and speech recognition.
For high-dimensional unstructured data for function regression, we focus in the forthcoming on plain feedforward (residual) neural networks, because neither a spatial structure (as in image data) nor a temporal structure (as in time series) is assumed.
Standard feedforward networks are given as a concatenation of $L$ consecutive layers, each described by a simple nonlinear transform:
    \begin{align}
    \label{eq:nn_one_layer_mapping}
        f_l(x) = \sigma(W_l x + b_l).
    \end{align}
As it turns out, the feedforward NN can suffer some numerical instabilities that can be partly elevated by the so-called residual network architecture where each layer takes the form
        \begin{align}
    \label{eq:nn_one_layer_mapping_resnet}
        f_l(x) = x+ \sigma(W_l x + b_l).
    \end{align}
The function $\sigma$ is called activation function, which is a pointwise acting nonlinear function that introduces nonlinearity into the approach.
A prominent example is the ReLU function $\sigma(x) = \max(x, 0)$ or variants like Leaky ReLU, SeLU, GeLU among others.
The matrices $W_l$ and the bias vectors $b_l$ constitute the parameters of the neural network, which can be optimized.
Due to the highly-nonlinear structure of the NN, a closed form solution for the optimal weights $\{ W_l \}_{l=1}^L$ and biases $\{ b_l \}_{l=1}^L$ is usually infeasible.
Thus, in practice neural networks are optimized using gradient descent like mini-batch optimization strategies, where the Adam optimizer \cite{kingma2014adam} which uses adaptive moment estimations is probably the most popular optimizer.
In order to improve generalization for neural networks, regularization strategies like DropOut and Early Stopping are available,
although modern neural network models are frequently overparameterized, mitigating the need for explicit regularization.
Despite there are error bounds available which elucidate especially the benefits of deep neural networks \cite{telgarsky2015representation}, 
they are rather of constructive nature.
As it is highly unlikely that the stochastic gradient based optimization realizes such constructions, these error bounds are usually not practical.
It is remarkable, that recently connections between the training of neural network and kernel methods where found: 
In the limit of large width neural networks, the training behaviour of neural networks is linearized in parameter space and can thus be described with help of the neural tangent kernel \cite{jacot2018neural}. 




\section{Numerical tests}\label{sec:NE}

In this section, we embark on an exploration and comparison of the supervised learning techniques introduced in the previous sections through a series of numerical tests.
In all the examples we employ sample-based algorithms for function regression. 
Given a set of sample points $\{x_j\}_{j=1}^n$, 
we measure the error between the prediction $s(x_j)$ and the true target value $y_j=f(x_j)$ considering the following relative error in norm 2:
\begin{align}
\label{eq:test_error}
err_{2} = \sqrt{\frac{\sum_{j=1}^n |y_j - s(x_j)|^2}{\sum_{j=1}^n |y_j|^2}}.
\end{align}
In the tables below, we add \textit{train} or \textit{test} in the subscript in order to define whether the error was measured with respect to the training or the test set.
In the last example we will consider also the error in the computation of the total cost obtained along the optimal trajectories. 
To this end, we introduce $cost_{SDRE}(x)$ as the total cost computed using the SDRE feedback and $cost_{surr}(x)$ as the total cost calculated 
using the gradient of the surrogate to the value functional as control signal in \eqref{optc}.
Finally, we define the error in the total cost at a given initial condition $x \in \R^d$:
$$
err_{cost}(x) = |cost_{SDRE}(x)-cost_{surr}(x)|.
$$

\begin{remark}
    Notice that the methods were running on different machines and with different languages. 
    Especially, the TT cross methods is based on a very performative implementation,
    while the block-sparse TT and kernel approximation did not use such tuned implementations.
    Moreover, the TT Gradient Cross and its applications are Matlab-based, while the rest of the codes are written in Python and Julia.
\end{remark}


\subsection{Low against high rank}

\label{low_rank_test}

First, we consider the following three test cases

\begin{enumerate}
\item[a)] $f(x)=\exp(-\sum_{i=1}^d x_i/(2d) )$, \quad  $x\in [-1,1]^d$,
\item[b)] $f(x)=\exp(-\prod_{i=1}^d x_i )$, \quad $x\in [-1,1]^d$,
\item[c)] $f(x)=\exp(-\prod_{i=1}^d x_i )$, \quad $x\in [0,2]^d$.
\end{enumerate}
We note that $f(x)=\exp(-\sum_{i=1}^d x_i/(2d) )$ has a rank one decomposition in a functional tensor train in contrast to the other two examples.  
The functions reported in the cases $(b)$ and $(c)$ are not low-rank and this is reflected in the numerical results. 
The difference between cases $(b)$ and $(c)$ are the different domains. 
Test case $(b)$ is set in $[-1,1]^d$ and the function is almost 1 for a significant part of the domain, while in case $(c)$ the function takes values between 0 and 1.

For each of three test cases, the methods described in Section \ref{sec:ML} were applied. 
The approximation take place in dimension $d=16$ and the tables for each methods display the train and test error as well as the degrees of freedom (DoFs), 
the CPU time and the sample size. 
The test error is based on a $10^4$ uniform random samples in the respective domain.

For the TT Gradient Cross, the space domain is discretized with 7 Legendre-Gauss nodes for each direction. We fix the stopping error $tol_{stop}$ equal to $10^{-5}$.
For the block sparse TT we consider degree 7 Legendre polynomials and a locality to match the degrees of freedom of the TT Cross results. 
Throughout the experiments, the stopping criterion was set to $10^{-11}$ or a maximal ALS iteration.
For the kernel approximation, we considered 5000 low discrepancy points within the given domain, and use the quadratic Matérn kernel with shape parameter in $\{ \frac{1}{2\sqrt{d}}, \frac{1}{4\sqrt{d}}, \frac{1}{8\sqrt{d}} \}$ and without explicit regularization. 
For the neural network, as standard setup and training was used: 
The neural networks used three fully connected layers of width 512,
and the training was done using mini batches of size 128, the Adam optimizer and an initial learning rate of 1e-3.

As the function considered in case $(a)$ has an exact TT decomposition with rank 1, it can be approximated arbitrary well by the TT Gradient Cross, 
as described by both $err_{train,2}$ and $err_{test,2}$. 
Also the low number of samples and the small CPU times are due to the low rank structure of the function. 
Similar, the block-sparse TT methods behaves very well. 
The cases $(b)$ and $(c)$ are more tricky and the magnitudes of the error indicators drops of several orders with respect to the first case. 
Especially for the third case for which the number of training samples increases drastically, 
obtaining just a testing error of order $O(10^{-1})$ for the block sparse approach and $O(10^{-2})$ for the TT Gradient Cross. 
For the block sparse TT format, an additional difficulty is due to the functions $(b)$ and $(c)$ not satisfying a locality bound. 
As no rounding was employed in the block sparse format, the number of degrees of freedom might be an overestimation. 
The kernel methods also are able to recover the chosen functions, despite the errors are not as small.
Especially the function $(c)$ is the most challenging one with an error in $O(10^{-1})$. 
This numbers are an indicator, that the kernel approximation using a radial kernel together with uniformly distributed samples is clearly affected from the curse of dimensionality.
Note that the kernel methods use significantly less samples and degrees of freedom. 
The neural networks achieve an error of $O(10^{-1})$ or $O(10^{-2})$, which is frequently sufficient for machine learning purposes, 
however here not on par with the kernel methods or the tensor trains.
Due to the stochastic iterative training procedure of the neural networks, 
some time had to be spent to find suitable hyperparameters, 
while there are likely still further options to tune the NN and its training.
Overall, the TT Cross approximation provides the best approximations.

\begin{table}[htbp]							
\centering
\begin{tabular}{c c  |  c c c c c c c}
\toprule 
&& $err_{train,2}$ & $err_{test,2}$ & DoFs & \# train sample & CPU train (s) & CPU test (s) \\
\midrule
TT Cross & (a) 
& 6.97e-16 & 7.60e-16 & 420 & 1253 &  0.11 &  0.47 \\
&(b) 
& 1.56e-7 & 3.16e-8 & 3612 & 13937 &  0.23 & 0.33 \\
&  (c) 
  & 7.66e-4 & 2.76e-2 & 35196 & 293265 &  3.45 & 0.41 \\

\midrule
BSTT &(a) & 6.90e-11 & 1.28e-09 & 434& 1900 & 27.5&0.15\\
&(b) & 4.40e-5 & 1.54e-4 &3689 & 14000& 1262.5 &0.60 \\
&  (c) & 6.03e-3 &9.05e-2 &27935 &100000&43973.7& 28.0\\

\midrule

NN&(a)  & 1.23e-2 & 1.56e-2 & 797185 & 5000 & 46.48 & 0.13 \\
&(b)    & 9.39e-3 & 1.41e-2 & 797185 & 5000 & 56.04 & 0.12 \\
&(c)    & 6.32e-2 & 1.14e-1 & 797185 & 5000 & 32.12 & 0.11 \\

\midrule
Kernel & (a)    & 4.76-12   & 9.35e-6  & 5000 & 5000 & 2.30 & 2.54 \\
 & (b)          & 2.98e-10  & 9.63e-4  & 5000 & 5000 & 2.39 & 2.60 \\
 & (c)          & 2.98e-10  & 3.61e-1   & 5000 & 5000 & 2.36 & 2.55 \\
\bottomrule		
\end{tabular}	
\caption{Results of the different methods for Test \ref{low_rank_test} $(d=16)$.}
\end{table}

\subsection{Regularity test}
\label{subsec:reg_test}

For the next test, we are interested in the effects of non-differentiablities of the target functions on the expressability of our model classes. To that end,  we consider the following family of functions
\begin{align*}
f(x) = \lambda_0 \Vert x \Vert^2 + \lambda_1 \Vert x -y_1 \Vert + \lambda_2 \sqrt{\Vert x -y_2 \Vert}, \quad x\in [-1,1]^d,    
\end{align*}
where $y_1 = (0.5,\ldots, 0.5)$ and $y_2 = (-0.5,\ldots, -0.5)$ and consider the cases of 
$\lambda = (\lambda_0,\lambda_1,\lambda_2) \in \{ (1,0,0), (1,0.5,0), (1,0.5,0.5), (0,0.5,0), (0,0,0.5)\}$.
Again, all the method described in Section \ref{sec:ML} are applied.

For the TT Gradient Cross, the intervall $[-1, 1]$ is discretized using 7 Legendre-Gauss nodes per direction and a stopping tolerance of $tol_{stop}= 10^{-5}$ is fixed. The block sparse TT where chosen to have the same complexity as the TT Cross approximation.
Both the kernel approximation and the neural network are using the same setup as in the previous example.
For every method, the test error is evaluated on $10^4$ random samples in $[-1,1]^d$.

\begin{figure}[htbp]\label{fig:diff}	
\centering
     \includegraphics[width=0.45\textwidth]{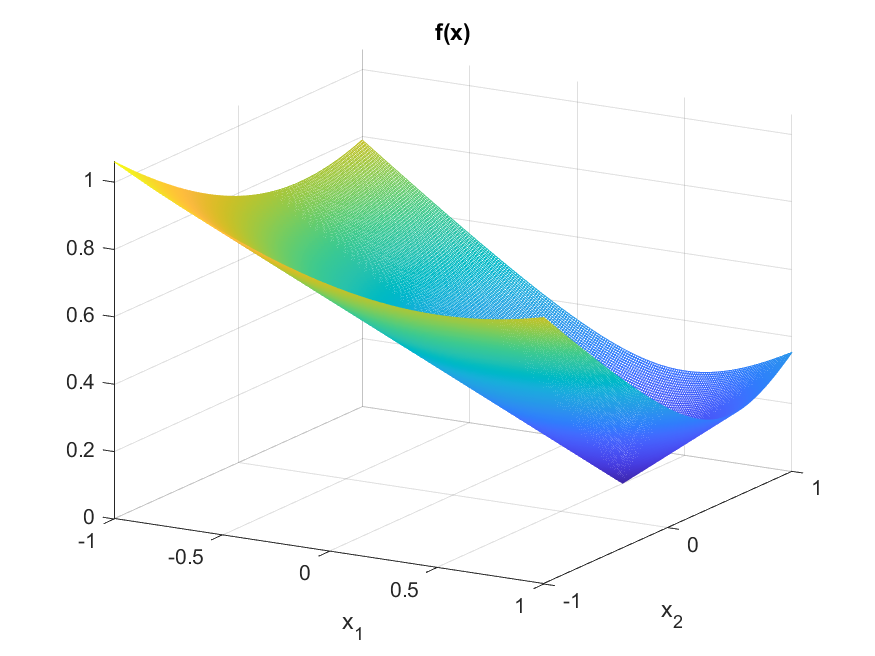}
    \includegraphics[width=0.45\textwidth]{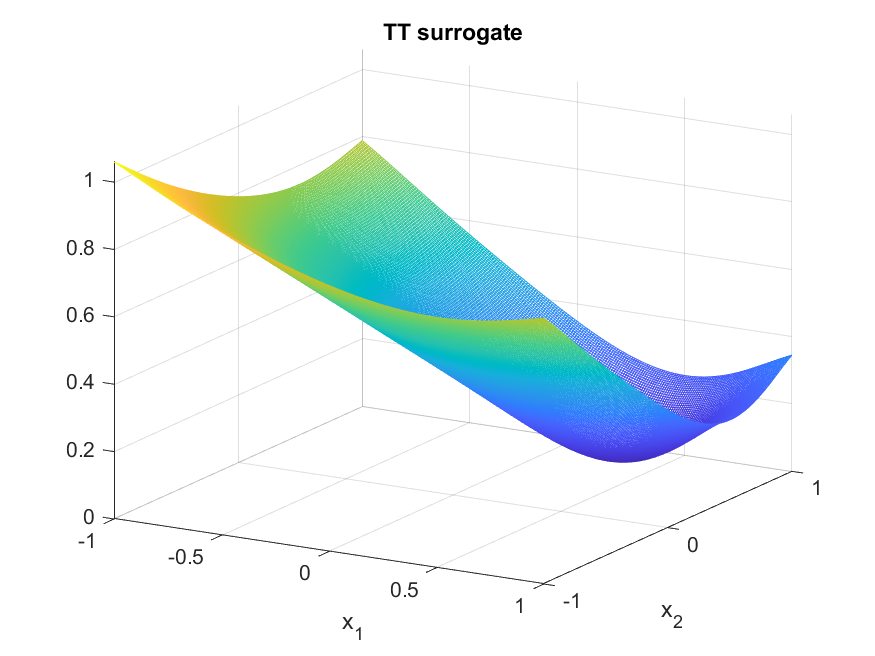}
     \includegraphics[width=0.45\textwidth]{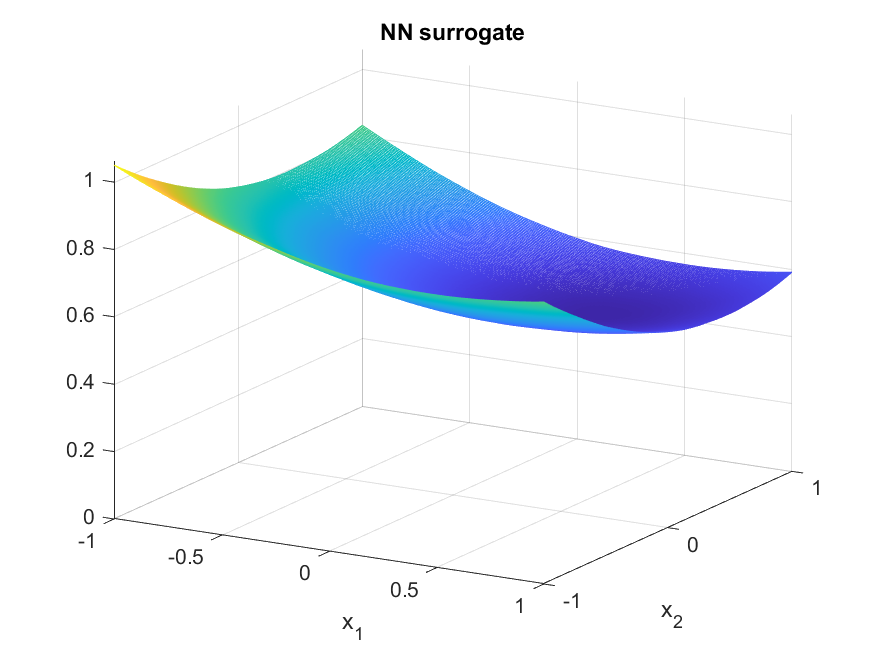}
 \includegraphics[width=0.45\textwidth]{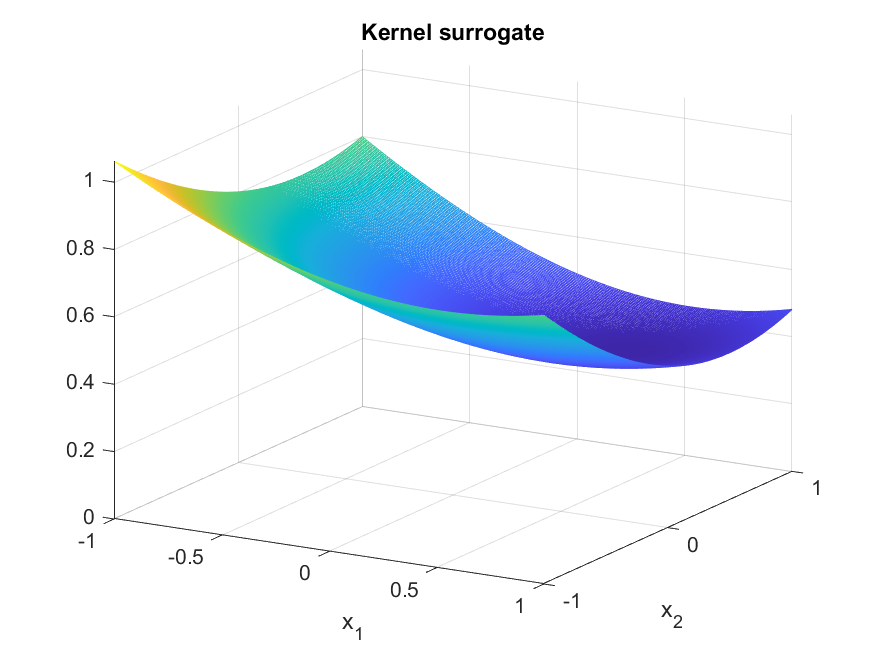}

	\caption{Function and surrogate model for $d=16$ and $\lambda = (0,0.5,0)$ for the different methods on the plane $(x_1,x_2,0.5,\ldots,0.5)$.
 The localized kink is only fitted properly by the TT approximation.}
\end{figure}

In this example one can see that, as expected, non-differentiablilities are more difficult to approximate with all of our methods.
Here, again the TT Cross method outperforms the other model classes providing results with improved accuracy of two orders of magnitude compared to the block sparse format and the kernel methods and up to four orders of magnitude compared to the NN approach. 
The TT methods can recover the squared norm almost exactly. 
Interestingly, there is almost no difference between the approximation error for $f(x)=\|x\|$ and $f(x)=\sqrt{\|x\|}$ within the different model classes respectively. 
This is likely due to the fact that the non-differentiability is very localized and thus not properly detected by the $10^4$ randomly sampled test points.
This is displayed in Figure \ref{fig:diff}. 
We note that TT-Cross algorithm is able to mimic the behaviour around the kink, while the kernel and NN approximation simply leave out this localized kink.

\begin{table}[htbp]							
		\centering
		\begin{tabular}{c|c c c|  c c c c c}
					\toprule
				 & $\lambda_0$ & $\lambda_1$ & $\lambda_2$
    & $err_{train,2}$ & $err_{test,2}$ & DoFs & \# train samples 
					     \\
					\midrule
    TT Cross&1 & 0 & 0 
    & 7.50e-16 & 9.64e-16 & 924 & 1994 \\
    &1 & 0.5 & 0 
    & 7.67e-7 & 1.57e-6 & 6608 & 14879 \\
    &1 & 0.5 & 0.5 
    & 1.49e-6 & 3.33e-6 & 13293 & 32235 \\
    &0 & 0.5 & 0 & 
    2.85e-7 & 2.19e-6 & 7322 & 16776 \\
    &0 & 0 & 0.5 
    & 5.83e-7 & 2.56e-6 & 7308 & 17334 \\

					\toprule

	BSTT&1 & 0 & 0
 & 9.62e-12 & 4.57e-11 & 96 & 1994 \\
    &1 & 0.5 & 0  & 1.41e-05  & 5.02e-5 & 5069 & 10000 \\
    &1 & 0.5 & 0.5   & 2.57e-05 &   7.4e-5  & 12786 & 32000 \\
   & 0 & 0.5 & 0  & 8.56e-05  & 3.4e-4  &  7378 & 16000 \\
   & 0 & 0 & 0.5  &  8.12e-05 & 3.1e-4   & 7378 & 17334 \\
 			\bottomrule		
			
NN&1 & 0 & 0 &        1.86e-2  & 3.13e-2 & 797185  & 5000 \\
&1 & 0.5 & 0 &      1.501-2 & 2.62e-2 & 797185 & 5000 \\
&1 & 0.5 & 0.5 &    9.82e-3 & 2.17e-2 & 797185 & 5000 \\
&0 & 0.5 & 0 &      4.78e-3 & 9.05e-3 & 797185 & 5000 \\
&0 & 0 & 0.5 &      8.41e-3 & 1.11e-2 & 797185 & 5000 \\
 			\bottomrule	
			
Kernel&1 & 0 & 0 &         9.74e-8 & 1.07e-4 & 5000 & 5000 \\
&1 & 0.5 & 0 &       9.01e-9 & 2.83e-4 & 5000 & 5000\\
&1 & 0.5 & 0.5 &     3.75e-7 & 2.18e-4 & 5000 & 5000\\
&0 & 0.5 & 0 &       4.34e-13 & 1.14e-3 & 5000 & 5000\\
&0 & 0 & 0.5 &       4.28e-13 & 1.08e-3 & 5000 & 5000\\
 			\bottomrule		
		\end{tabular}	
\caption{Results for the different methods for \Cref{subsec:reg_test} $(d=16)$.}
\end{table}

\subsection{Academic Optimal Control Example}

In this numerical example we want to deal with the application of the mentioned techniques to an ``academic'' optimal control problem. 
Following the idea described in \cite{ehring2023hermite}, our aim is to consider a control problem with a prefixed value function $V$.
Now, consider the dynamical system 
\begin{equation}
    \dot x = u, \quad x(0)=x_0
    \label{dyn_academic}
\end{equation}
 and the cost functional 
 \begin{equation}
 J(x_0,u) = \int_0^{+\infty} r(x(t)) + \frac{1}{2}\|u(t)\|^2\mathrm{d}t,
 \label{cost_academic}
 \end{equation}
where $r(x)$ will be chosen to obtain the value function $V$.
In this case the HJB equation reads
\begin{align*}
-\frac{1}{2}\|\nabla V(x)\|^2 + r(x)= 0.
\end{align*}
Then, choosing exactly $r(x) = \frac{1}{2}\|\nabla V(x)\|^2$, 
we obtain that the optimal control problem \eqref{dyn_academic}-\eqref{cost_academic} admits $V$ as value function. 
The aim of this example is twofold: 
We want to test the accuracy of the SDRE approach in the approximation of the feedback law and we want to compare the supervised learning techniques varying the dimension of the problem.

\paragraph{SDRE approximation}

In this paragraph we focus on the accuracy of the SDRE approximation of the optimal control problem \eqref{dyn_academic}-\eqref{cost_academic}.
The dynamical system \eqref{dyn_academic} can be written easily in semilinear form \eqref{semilinear} with $A(x) = 0_d$ and $B(x) = I_d$, where $0_d \in \R^{d \times d}$ is a matrix of all zeros and $I_d \in \R^{d \times d}$ is the identity matrix. The cost functional \eqref{cost_academic} is then written in a state-dependent quadratic form \eqref{quadratic_cost} with $R(x) = \frac{1}{2}I_d$ and fixing $Q(x)$ such that
$$x^\top Q(x)x=r(x) = \frac 1 2\|V(x)\|^2 \; .$$
In this case the associated SDRE for $x \in \Omega$ reads
\begin{align*}
-2P(x)^2 + Q(x) = 0
\end{align*}
with solutions $ \pm \sqrt{Q(x)/2}$. 
Since we are interested in the positive definite solution, we choose $P(x) = \sqrt{Q(x)/2}$. 
We note that the choice of $Q(x)$ is crucial (see \cite{dolgov2022optimizing,Astolfi2020} for a discussion on the importance of the representation of the system). 
Indeed, choosing 
$$
  Q(x) = \frac{1}{2}
  \begin{bmatrix}
    |\partial_{x_1} V|^2/{x^2_1} & & \\
    & \ddots & \\
    & & |\partial_{x_d} V|^2/{x^2_d}
  \end{bmatrix} \, ,
$$ the solution of the SDRE reads
\begin{align*}
  P(x) = \frac{1}{2}
  \begin{bmatrix}
    \partial_{x_1} V/{x_1} & & \\
    & \ddots & \\
    & & \partial_{x_d} V/{x_d}
  \end{bmatrix} \; ,
\end{align*}
and the corresponding feedback control can be computed via the formula \eqref{control_SDRE}: 
$$
u(x) = -\frac{1}{2} P(x) x = -\frac{1}{4} \nabla V,
$$
retrieving exactly the formula of the optimal feedback control \eqref{optc}. 
In this specific example the exactness of the SDRE approach comes only in case we are able to extract the norms of the single partial derivatives $|\partial_{x_i} V|^2$ from the the given cost function $r(x) = \frac 1 2\|V\|^2$.

\paragraph{Supervised learning approximation}
Now let us compare the different techniques on a specific example. Let us choose the value function of the form
\begin{align}
    V(x) = \|x\|^2 (e^{-\frac{\|x-\mu_1\|^2}{\sigma_1^2}}+ e^{-\frac{\|x-\mu_2\|^2}{\sigma_2^2}}),
    \label{eq:vf_academic}
\end{align}
where $\mu_1,\mu_2 \in \mathbb{R}^d$ and $\sigma_1,\sigma_2 \in \mathbb{R}^+$ are parameters denoting respectively the means and the standard deviations of the corresponding gaussian functions. We can see immediately that the value function can be expressed in the form $V(x) = x^\top P(x) x$ defining
$$
  P(x) = 
  \begin{bmatrix}
  (e^{-\frac{\|x-\mu_1\|^2}{\sigma_1^2}}+ e^{-\frac{\|x-\mu_2\|^2}{\sigma_2^2}}) & & \\
    & \ddots & \\
    & & (e^{-\frac{\|x-\mu_1\|^2}{\sigma_1^2}}+ e^{-\frac{\|x-\mu_2\|^2}{\sigma_2^2}})
  \end{bmatrix} \; .
$$

First of all, we fix $\mu_1 = 0 \cdot \mathbf{1}$, $\mu_2 = 0.5 \cdot \mathbf{1}$, $\sigma_1 = \sigma_2 = 1$, where $\mathbf{1} \in \mathbb{R}^d$ is a vector with all ones. 
The dimension of the problem $d$ will vary in the set $$D = \{3,\ldots,16\}$$ and the initial condition $x_0$ will be taken randomly in the set $\Omega = [-1,1]^d$. 
As the variance is dimension independent one would expect the function to become smoother with increased dimension. 
Indeed, the TT Cross method has an improving approximation error with respect to the dimensions.

The TT-Gradient Cross is implemented using a stopping tolerance $10^{-5}$ and 7 Legendre-Gauss polynomials per dimension. 
We see that as we increase the dimension the error improves, but stagnating around the order $10^{-3}$. 
On the hand, the CPU is keeping an almost constant behaviour of order $O(1)$, due to the fact that the TT rank is at most 5 for all the considered dimensions.


For the kernel approximation, we employ the quadratic Matérn kernel using 10000 low-discrepancy points within $[-1, 1]^d$ and choose the shape parameter as $\frac{1}{2\sqrt{d}}$. 
The kernel approximation is the best approximation up to dimension 5, while its error increases as the dimension grows. 
The behaviour of the CPU time is more or less constant of $O(10)$, as it mainly depends on the number of sample points. 
Finally, the NN approach, which uses the same setup as in the previous examples, obtains the worst results in terms of both indicators, with an increasing error between $O(10^{-2})$ and $O(10^{-1})$ as the dimension $d$ of the problem grows and a CPU time varying between $O(10)$ and $O(10^2)$.

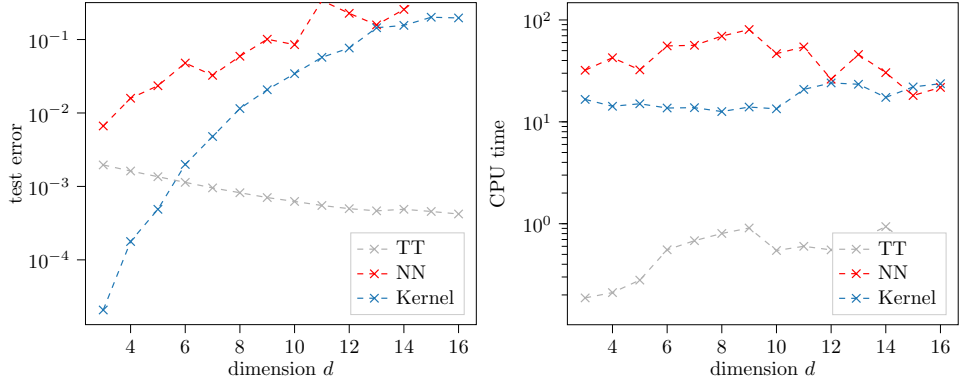
\begin{figure}[h]
\centering
\begin{minipage}[t]{.475\textwidth}
\begin{tikzpicture}[scale=0.75]

\definecolor{darkgray176}{RGB}{176,176,176}
\definecolor{lightgray204}{RGB}{204,204,204}
\definecolor{steelblue31119180}{RGB}{31,119,180}

\begin{axis}[
legend cell align={left},
legend style={
  fill opacity=0.8,
  draw opacity=1,
  text opacity=1,
  at={(0.97,0.03)},
  anchor=south east,
  draw=lightgray204
},
log basis y={10},
tick align=outside,
tick pos=left,
x grid style={darkgray176},
xlabel={dimension $d$},
xmin=2.35, xmax=16.65,
xtick style={color=black},
y grid style={darkgray176},
ymin=1.31247271787062e-05, ymax=0.318279102007188,
ylabel={test error},
ymode=log,
ytick style={color=black},
ytick={1e-06,1e-05,0.0001,0.001,0.01,0.1,1,10},
yticklabels={
  \(\displaystyle {10^{-6}}\),
  \(\displaystyle {10^{-5}}\),
  \(\displaystyle {10^{-4}}\),
  \(\displaystyle {10^{-3}}\),
  \(\displaystyle {10^{-2}}\),
  \(\displaystyle {10^{-1}}\),
  \(\displaystyle {10^{0}}\),
  \(\displaystyle {10^{1}}\)
}
]
\addplot [semithick, darkgray176, dashed, mark=x, mark size=3, mark options={solid}]
table {%
3   0.001961249937931
4   0.001621245829460
5   0.001356435306973
6   0.001135330918923
7   0.000954612208148
8   0.000820521324684
9   0.000705082624050
10  0.000624319516931
11  0.000550161027489
12  0.000497922017728
13  0.000465043603165
14  0.000488146527062
15  0.000455756738550
16  0.000421249659030
};
\addlegendentry{TT}
\addplot [semithick, red, dashed, mark=x, mark size=3, mark options={solid}]
table {%
3   0.006654322865737672
4   0.015865769158239524
5   0.0234762699647523
6   0.04788314698413772
7   0.03251551831233177
8   0.059205290779729705
9   0.10096277249621967
10  0.08511370055561222
11  0.33930733083866793
12  0.22737139856828736 
13  0.15818115727977197
14  0.25696650611137595
15  0.49767833929785865
16  0.6194948926575076
};
\addlegendentry{NN}
\addplot [semithick, steelblue31119180, dashed, mark=x, mark size=3, mark options={solid}]
table {%
3 2.07680461843475e-05
4 0.000177768783155974
5 0.000487918248096482
6 0.00199251839930451
7 0.00477994955953652
8 0.0114851514194671
9 0.0207228748704445
10 0.0342923525426094
11 0.057236650968847
12 0.0765408324450994
13 0.144794175680928
14 0.155816101108935
15 0.201142001681233
16 0.196798993458981
};
\addlegendentry{Kernel}
\end{axis}

\end{tikzpicture}
\end{minipage}
\begin{minipage}[t]{.475\textwidth}
\begin{tikzpicture}[scale=0.75]

\definecolor{darkgray176}{RGB}{176,176,176}
\definecolor{lightgray204}{RGB}{204,204,204}
\definecolor{steelblue31119180}{RGB}{31,119,180}

\begin{axis}[
legend cell align={left},
legend style={
  fill opacity=0.8,
  draw opacity=1,
  text opacity=1,
  at={(0.97,0.03)},
  anchor=south east,
  draw=lightgray204
},
log basis y={10},
tick align=outside,
tick pos=left,
x grid style={darkgray176},
xlabel={dimension $d$},
xmin=2.35, xmax=16.65,
xtick style={color=black},
y grid style={darkgray176},
ylabel={CPU time},
ymode=log,
ytick style={color=black},
]
\addplot [semithick, steelblue31119180, dashed, mark=x, mark size=3, mark options={solid}]
table {%

};
\addplot [semithick, darkgray176, dashed, mark=x, mark size=3, mark options={solid}]
table {%
3   0.187024000000000
4   0.210797800000000
5   0.279492800000000
6   0.554857800000000
7   0.681390400000000
8   0.802043600000000
9   0.907791500000000
10  0.545695900000000
11  0.599964400000000
12  0.553491100000000
13  0.694471300000000
14  0.937298500000000
15  0.579127300000000
16  0.489735400000000
};
 \addlegendentry{TT}
\addplot [semithick, red, dashed, mark=x, mark size=3, mark options={solid}]
table {%
3   32.07794642448425
4   42.60289764404297
5   32.357810497283936
6   55.62722897529602
7   56.391977310180664
8   69.41751670837402
9   80.51443552970886
10  46.71954417228699
11  54.40788173675537 
12  26.210535526275635
13  45.80554676055908
14  30.40624761581421
15  18.132041215896606
16  21.76535391807556
};
\addlegendentry{NN}
\addplot [semithick, steelblue31119180, dashed, mark=x, mark size=3, mark options={solid}]
table {%
3   16.600989818572998
4   14.21419095993042
5   15.04184365272522
6   13.66795539855957
7   13.778729677200317
8   12.615983009338379
9   13.970061302185059
10  13.419519424438477
11  20.821945190429688
12  24.089510679244995
13  23.376524925231934
14  17.364453315734863
15  22.024975299835205
16  23.67611789703369
};
\addlegendentry{Kernel}
\end{axis}

\end{tikzpicture} 
\end{minipage}
\caption{Visualization of the test error \eqref{eq:test_error} (left) and the CPU time (right) for approximation of the value function \eqref{eq:vf_academic} depending on the dimension $d$ ($x$-axis).}
\label{fig:academic_example}
\end{figure}

\subsection{Control of the Allen-Cahn equation}
Let us consider the following nonlinear Allen-Cahn PDE with homogeneous Neumann boundary conditions:
\begin{equation*}
\left\{ \begin{array}{l}
\partial_t y(t,x) = \sigma \partial_{xx} y(t,x) + y(t,x) (1-y(t,x)^2) + u(t,x),  \\
 y(0,x)=y_0(x),
\end{array} \right.
\label{AC}
\end{equation*}
with $x \in [0,1]$ and $t \in (0,+\infty)$ and the following cost functional
$$
J(u,y_0) = \int_0^{\infty}  \int_0^1 (|y(t,x)|^2 +  \gamma |u(t,x)|^2) dx \, dt \,.
$$
Approximating the PDE by finite difference schemes with $d$ grid points, we obtain the following ODEs system 
\begin{align*}
\dot{y}(s)=A(y) y(s) + u(s),
\end{align*}
with
\begin{align*}
A(y) = \sigma A_0 + I_d -  diag(y \odot y),\quad y \in \mathbb{R}^d, 
\end{align*}
where $\odot$ is the Hadamard product, $I_d \in \mathbb{R}^{d \times d}$ is the identity matrix and $A_0$ is the tridiagonal matrix arising from the discretization of the Laplacian with Neumann boundary conditions. We fix $\sigma = 10^{-2}$ and  $d=30$. The different techniques will be tested on a set of initial conditions of the form
\begin{equation}
y_0(x) = \sum_{k=1}^{4} \frac{a_k}{2} \cos(2 \pi k x) k^{-\beta}, \quad a_k \in \{0,1\},
\label{y0_fourier}
\end{equation}
where the parameter $\beta$ is related to the
decay of the Fourier coefficients and to the regularity we want to assume. We fix $\beta=3$.
The vector $x=[x_1,\ldots, x_d]$ contains the discretization points of the interval $[0,1]$.

The supervised learning techniques will be trained over the domain $[-1,1]^d$, while for the test phase we are going to consider 4 initial conditions in the form \eqref{y0_fourier} varying the vector $a=[a_1, \, a_2, \, a_3, \, a_4]$. To evaluate the feedback law induced by the surrogates of the test functions we integrate the dynamical systems with these controls until a final time $t_{\mathrm{final}} =60$ and use the trapezoidal rule to estimate the costs along these trajectories. We compare this to the costs induced by the SDRE feedback in $err_{cost}$.


\vspace{5mm}

We start applying the TT Gradient Cross for the approximation of the value function. We fix $n=6$ Gauss-Legendre nodes per dimension and a stopping tolerance $tol_{stop} = 10^{-4}$. After 17 iterations the algorithm stopped reaching the prescribed tolerance and the final TT rank $r$ is equal to 18. Table \ref{table_train_AC} reports the errors, the computational cost and the number of samples used during the training phase. We note that for the error $err_{train,2}$ for TT Gradient Cross is below the prescribed tolerance ($10^{-4}$) and it completes the training of the surrogate model in just almost 20 seconds. The number of training samples is comparable to the complexity of the TT structure which we recall to be equal to $O(dnr^2)$.

The kernel approximation uses the Gaussian kernel with shape parameter $1/\sqrt{d}$, which gave slightly better results than the use of the quadratic Matérn kernel.
The neural network uses the same setup and training hyperparameters as in the previous examples.
Both the kernel approximation and the neural network are trained using 5000 samples, which were generated randomly in Fourier space.
The kernel approximation accuracy is comparable to the TT approximation, while the NN performs worse.
In terms of the training time, the kernel methods performs best, while the NN has a long training time due to its iterative stochastic gradient descent optimization.
We note that the kernel approximation performs comparably well here compared to the TT, because it is only trained on meaningful inputs which were sampled in Fourier space, while the TT likely provides a suitable approximation of the value function in the whole domain.


\begin{table}[hbht]
\centering
\begin{tabular}{c|c|c|c}    
Surrogate model & $err_{train,2}$ & CPU train ($s$) & \# train samples \\ 
\hline
TT Gradient Cross &    3.48e-5 &  20.10 & 121906 \\
NN & 2.03e-2 & 45.80 & 5000 \\
Kernel & 2.147e-5 & 2.82 & 5000 \\

 \end{tabular}
  \caption{Comparison in terms of the training phase for the different supervised learning techniques.}
 \label{table_train_AC}
\end{table}

Table \ref{table_test_AC} compares the different surrogate models in the testing phase using 4 initial conditions in the form \eqref{y0_fourier} varying the vector $a \in \R^4$. 
As regards TT Gradient Cross, 
the $err_{test}$ reaches the same order of magnitude of the prescribed stopping tolerance for all the studied cases, 
while for the $err_{cost}$ we loose two orders of magnitude with respect to the $err_{test}$, due to the fact that the error in the cost takes into account the construction of the control, hence depending on the gradient of the surrogate models. 
The kernel approximation and the neural network fail to stabilize the dynamical system in a neighbourhood of the origin:
More specifically, when the norm of the trajectories is below a certain threshold denoted by $a_{TB}$, the error in the approximation of the gradient leads to the synthesis of a non-stabilizing feedback control and the solution is driven far away from the origin. 
To this end, we consider a modification in the computation of the feedback control already introduced in \cite{dolgov2022data} and denoted as Two-Boxes (TB) approach. 
More precisely, in the region in which the surrogate is not able to stabilize the system, 
we substitute the surrogate control with a control given by the Linear Quadratic Regulator, 
i.e. the solution of the SDRE at the origin. 
Denoted by $P_0$ the solution of the Riccati equation corresponding to the Linear Quadratic Regulator, we define the surrogate-TB feedback control as:
\begin{align*}
u^*(x)= \left\{ \begin{array}{l}
-R^{-1} B^{\top} Px, \qquad \qquad \qquad \Vert x \Vert \le  a_{TB} ,\\
-\frac{1}{2}R^{-1} B(x)^{\top} \nabla V_{surr}(x), \quad \Vert x \Vert >  a_{TB} ,
\end{array} \right.
\end{align*}
where $\nabla V_{surr}(x)$ is the gradient of the surrogate model for the approximation of the value function.

Using this strategy, we note by Table \ref{table_test_AC} that TT and Kernel-TB achieve almost the same $err_{cost}$, 
while NN-TB gets slightly higher results. 
Regarding the $err_{test}$ we observe that Kernel method achieves the best accuracy with order $O(10^{-6})$, while the NN approach gets better results increasing the number of the terms in the Fourier expansion \eqref{y0_fourier}.

\begin{table}[hbht]
\centering
\begin{tabular}{c|cc|cc|cc}     
  & \multicolumn{2}{c|}{TT Gradient Cross}    & \multicolumn{2}{c}{NN-TB} &  
\multicolumn{2}{|c}{Kernel-TB} \\ 
     & $err_{test}$ & $err_{cost}$ & 
     $err_{test}$ & $err_{cost}$ & $err_{test}$ & $err_{cost}$
\\\hline
$[1,0,0,0]$ & 1.77e-4  &  0.0320   & 1.90e-3 & 0.0545 & 1.32e-6 & 0.0341 \\
$[1,1,0,0]$ &   2.00e-4  & 0.0330  & 2.13e-3 & 0.0572 & 8.62e-6 & 0.0350 \\
$[1,1,1,0]$ & 2.00e-4   &  0.0333   & 8.11e-4 & 0.0577 & 1.24e-5 & 0.0360 \\
$[1,1,1,1]$ &   1.98e-4  &  0.0333  & 8.37e-5 & 0.0580 & 1.64e-5 & 0.0364 \\

 \end{tabular}
  \caption{Comparison in terms of the testing phase with 4 initial conditions in the form \eqref{y0_fourier} for the different supervised learning techniques.}
 \label{table_test_AC}
\end{table}





\begin{figure}[htbp]	
\centering
	\includegraphics[width=0.49\textwidth]{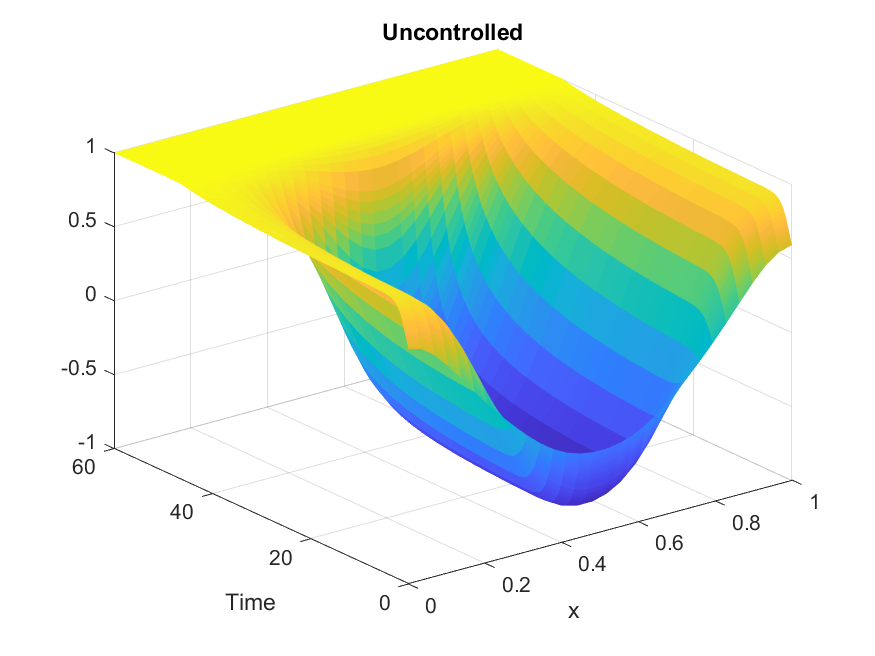}
 \includegraphics[width=0.49\textwidth]{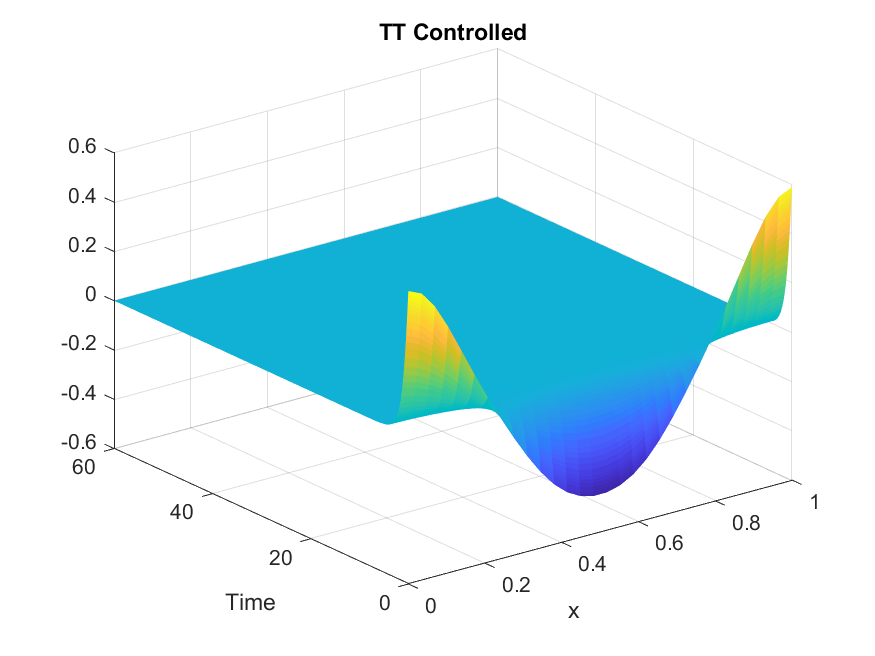}
\includegraphics[width=0.49\textwidth]{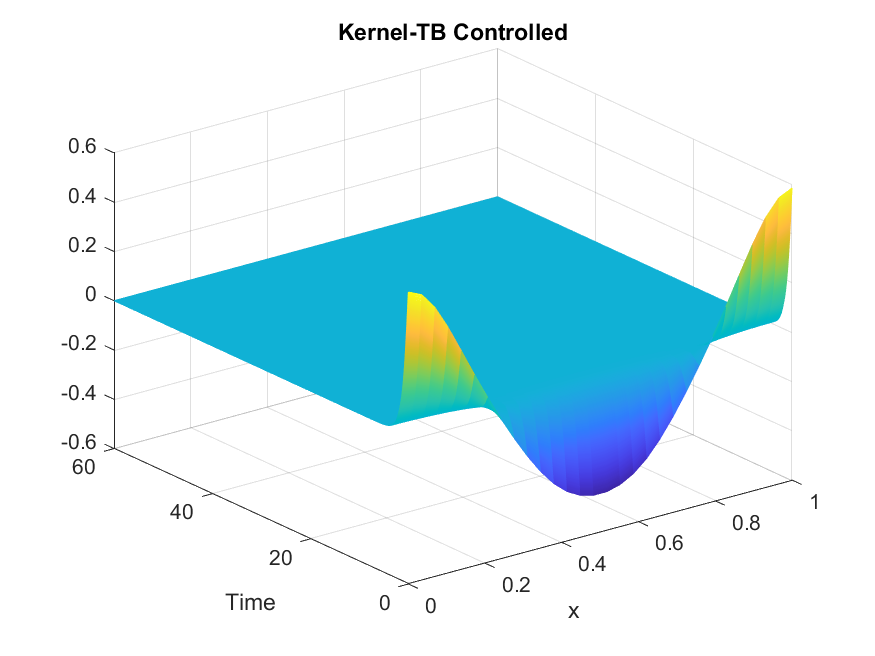}
 \includegraphics[width=0.49\textwidth]{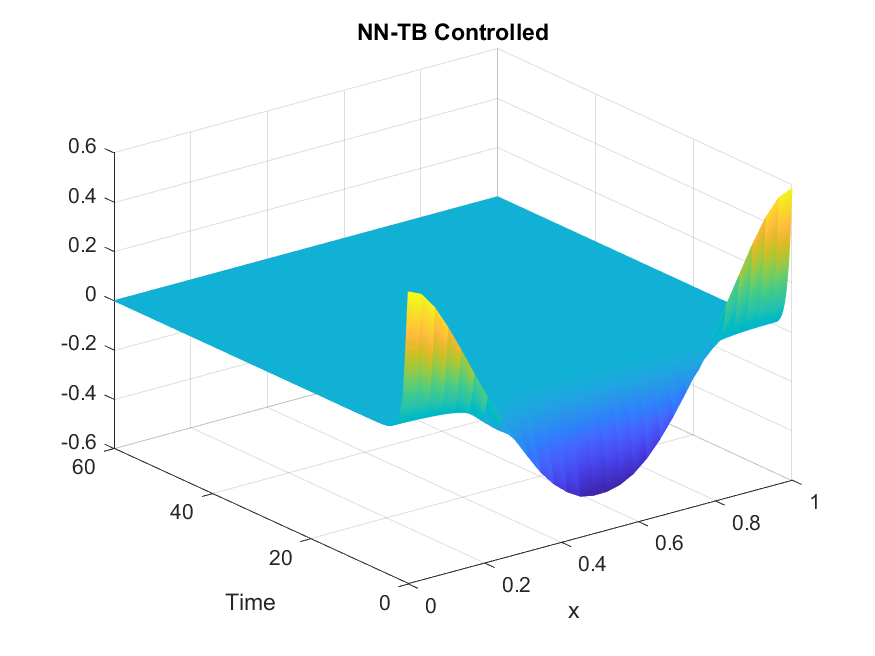}

	\caption{Solution for the Allen-Cahn equation for the uncontrolled dynamics (top left), the TT Gradient Cross controlled dynamics (top right), the Kernel-Two Boxes controlled dynamics (bottom left) and NN controlled dynamics (bottom right) }
 \label{AC_unc_TT}
\end{figure}

Finally, we show the plots of the uncontrolled and of the controlled solutions of the Allen-Cahn PDE.
The top left panel of Figure \ref{AC_unc_TT} displays the uncontrolled solution and we can note that at the final time it converges to the stable solution $\overline{y}_1(x) \equiv 1$. 
In the top right panel and in both lower panels of Figure \ref{AC_unc_TT} we can observe the behaviour of the Allen-Cahn solution controlled via the TT Gradient Cross, 
the Kernel-Two Boxes and the NN-Two Boxes approaches. 
We note that the control is driving the solution to the unstable equilibrium $\overline{y}_2(x) \equiv 0$ already at the very first time instances, 
keeping the solution close to $\overline{y}_2(x)$ for the remaining time instances.

\section{Conclusions}
\label{sec:conclusion}

Overall, supervised learning schemes based on tensor trains, kernel methods and neural networks are able to recover high dimensional functions in an error regime between $O(10^{-6})$ and $O(10^{-3})$ quite reliable. 
The experiments very nicely show the influence of intrinsic structures as regularity and separability (as used in functional tensor train ranks) on the approximation errors. 
As most of our test cases exhibit some kind of regularity, the TT methods performed very well, as they exploit such assumptions. 
The kernels models and neural networks could not quite beat tensor trains but allow for a much more general usage:
an important advantage of kernel methods and neural networks compared to the TT Cross is especially the ability to use scattered data for the approximation, while TT Cross depends on an active learning mode by adding samples based on a subset of grid points, which makes it difficult to employ on real world data.
One can also observe that the metric used to evaluate the functions plays an important role in the assessment of the methods, as the mean squared error used here does not necessarily identify non-differentiabilities or localizations. 

The use of available a priori information should be in the centre of high dimensional approximation tasks, when a reasonable high accuracy is required.

\section*{Acknowledgements}

This manuscript is dedicated to the memory of Maurizio Falcone. 
Maurizio was not only a great mathematician and an academic father, but also a father figure following all the steps of his academic "sons". 
He was really interested in high-dimensional optimal control problems and I hope he will enjoy this paper, whenever he is. 
Thank you for sharing all your knowledge, your passion for research and your smiles.

\section*{Declarations}

{\bf Ethical Approval} Not applicable.

\vspace{3mm}
\noindent
{\bf Competing interests} The authors have no competing interests.

\vspace{3mm}
\noindent
{\bf Authors' contributions} All the authors equally contributed to the preparation of the submitted manuscript.

\vspace{3mm}
\noindent
{\bf Funding} L. Saluzzi was supported by "Gruppo Nazionale per il Calcolo Scientifico" (GNCS - INdAM) and was "titolare di borsa per l’estero dell’Istituto Nazionale di Alta Matematica". M. Oster acknowledges funding by the Deutsche Forschungsgemeinschaft (DFG, German Research Foundation) - Project number 442047500 through the Collaborative Research Center "Sparsity and Singular Structures" (SFB 1481).

\vspace{3mm}
\noindent
{\bf Data Availability Statement} A Github repository will be made available upon publication of this work. No data has been used in this paper.

\bibliography{references} 


\begin{thebibliography}{93}
\ifx \bisbn   \undefined \def \bisbn  #1{ISBN #1}\fi
\ifx \binits  \undefined \def \binits#1{#1}\fi
\ifx \bauthor  \undefined \def \bauthor#1{#1}\fi
\ifx \batitle  \undefined \def \batitle#1{#1}\fi
\ifx \bjtitle  \undefined \def \bjtitle#1{#1}\fi
\ifx \bvolume  \undefined \def \bvolume#1{\textbf{#1}}\fi
\ifx \byear  \undefined \def \byear#1{#1}\fi
\ifx \bissue  \undefined \def \bissue#1{#1}\fi
\ifx \bfpage  \undefined \def \bfpage#1{#1}\fi
\ifx \blpage  \undefined \def \blpage #1{#1}\fi
\ifx \burl  \undefined \def \burl#1{\textsf{#1}}\fi
\ifx \doiurl  \undefined \def \doiurl#1{\url{https://doi.org/#1}}\fi
\ifx \betal  \undefined \def \betal{\textit{et al.}}\fi
\ifx \binstitute  \undefined \def \binstitute#1{#1}\fi
\ifx \binstitutionaled  \undefined \def \binstitutionaled#1{#1}\fi
\ifx \bctitle  \undefined \def \bctitle#1{#1}\fi
\ifx \beditor  \undefined \def \beditor#1{#1}\fi
\ifx \bpublisher  \undefined \def \bpublisher#1{#1}\fi
\ifx \bbtitle  \undefined \def \bbtitle#1{#1}\fi
\ifx \bedition  \undefined \def \bedition#1{#1}\fi
\ifx \bseriesno  \undefined \def \bseriesno#1{#1}\fi
\ifx \blocation  \undefined \def \blocation#1{#1}\fi
\ifx \bsertitle  \undefined \def \bsertitle#1{#1}\fi
\ifx \bsnm \undefined \def \bsnm#1{#1}\fi
\ifx \bsuffix \undefined \def \bsuffix#1{#1}\fi
\ifx \bparticle \undefined \def \bparticle#1{#1}\fi
\ifx \barticle \undefined \def \barticle#1{#1}\fi
\bibcommenthead
\ifx \bconfdate \undefined \def \bconfdate #1{#1}\fi
\ifx \botherref \undefined \def \botherref #1{#1}\fi
\ifx \url \undefined \def \url#1{\textsf{#1}}\fi
\ifx \bchapter \undefined \def \bchapter#1{#1}\fi
\ifx \bbook \undefined \def \bbook#1{#1}\fi
\ifx \bcomment \undefined \def \bcomment#1{#1}\fi
\ifx \oauthor \undefined \def \oauthor#1{#1}\fi
\ifx \citeauthoryear \undefined \def \citeauthoryear#1{#1}\fi
\ifx \endbibitem  \undefined \def \endbibitem {}\fi
\ifx \bconflocation  \undefined \def \bconflocation#1{#1}\fi
\ifx \arxivurl  \undefined \def \arxivurl#1{\textsf{#1}}\fi
\csname PreBibitemsHook\endcsname

\bibitem[\protect\citeauthoryear{Bellman}{1966}]{bellman1966dynamic}
\begin{barticle}
\bauthor{\bsnm{Bellman}, \binits{R.}}:
\batitle{Dynamic programming}.
\bjtitle{Science}
\bvolume{153}(\bissue{3731}),
\bfpage{34}--\blpage{37}
(\byear{1966})
\end{barticle}
\endbibitem

\bibitem[\protect\citeauthoryear{Bardi and
  Capuzzo-Dolcetta}{1997}]{DeterministicHJB}
\begin{bbook}
\bauthor{\bsnm{Bardi}, \binits{M.}},
\bauthor{\bsnm{Capuzzo-Dolcetta}, \binits{I.}}:
\bbtitle{Optimal Control and Viscosity Solutions of Hamilton-Jacobi-Bellman
  Equations}.
\bpublisher{Birk\"auser},
\blocation{Boston}
(\byear{1997})
\end{bbook}
\endbibitem

\bibitem[\protect\citeauthoryear{Falcone and Ferretti}{2013}]{falcone2013semi}
\begin{bbook}
\bauthor{\bsnm{Falcone}, \binits{M.}},
\bauthor{\bsnm{Ferretti}, \binits{R.}}:
\bbtitle{Semi-Lagrangian Approximation Schemes for Linear and Hamilton—Jacobi
  Equations}.
\bpublisher{SIAM},
\blocation{Philadelphia, PA}
(\byear{2013})
\end{bbook}
\endbibitem

\bibitem[\protect\citeauthoryear{Kalise and Kunisch}{2018}]{pol_approx_kunisch}
\begin{barticle}
\bauthor{\bsnm{Kalise}, \binits{D.}},
\bauthor{\bsnm{Kunisch}, \binits{K.}}:
\batitle{Polynomial approximation of high-dimensional
  {H}amilton-{J}acobi-{B}ellman equations and applications to feedback control
  of semilinear parabolic {PDE}s}.
\bjtitle{SIAM J. Sci. Comput.}
\bvolume{40}(\bissue{2}),
\bfpage{629}--\blpage{652}
(\byear{2018})
\end{barticle}
\endbibitem

\bibitem[\protect\citeauthoryear{Alla et~al.}{2015}]{alla2015efficient}
\begin{barticle}
\bauthor{\bsnm{Alla}, \binits{A.}},
\bauthor{\bsnm{Falcone}, \binits{M.}},
\bauthor{\bsnm{Kalise}, \binits{D.}}:
\batitle{An efficient policy iteration algorithm for dynamic programming
  equations}.
\bjtitle{SIAM Journal on Scientific Computing}
\bvolume{37}(\bissue{1}),
\bfpage{181}--\blpage{200}
(\byear{2015})
\end{barticle}
\endbibitem

\bibitem[\protect\citeauthoryear{Zhao et~al.}{2014}]{Zhao}
\begin{bchapter}
\bauthor{\bsnm{Zhao}, \binits{Z.}},
\bauthor{\bsnm{Yang}, \binits{Y.}},
\bauthor{\bsnm{Li}, \binits{H.}},
\bauthor{\bsnm{Liu}, \binits{D.}}:
\bctitle{Approximate finite-horizon optimal control with policy iteration}.
In: \bbtitle{Proceedings of the 33rd Chinese Control Conference},
pp. \bfpage{8895}--\blpage{8900}
(\byear{2014})
\end{bchapter}
\endbibitem

\bibitem[\protect\citeauthoryear{Tahirovic and Astolfi}{2019}]{linearlikePI}
\begin{bchapter}
\bauthor{\bsnm{Tahirovic}, \binits{A.}},
\bauthor{\bsnm{Astolfi}, \binits{A.}}:
\bctitle{Optimal control for continuous- time nonlinear systems based on a
  linear-like policy iteration}.
In: \bbtitle{2019 IEEE 58th Conference on Decision and Control (CDC)},
pp. \bfpage{5238}--\blpage{5243}
(\byear{2019})
\end{bchapter}
\endbibitem

\bibitem[\protect\citeauthoryear{He et~al.}{2020}]{OnlinePolicy}
\begin{barticle}
\bauthor{\bsnm{He}, \binits{S.}},
\bauthor{\bsnm{Fang}, \binits{H.}},
\bauthor{\bsnm{Zhang}, \binits{M.}},
\bauthor{\bsnm{Liu}, \binits{F.}},
\bauthor{\bsnm{Ding}, \binits{Z.}}:
\batitle{Adaptive optimal control for a class of nonlinear systems: The online
  policy iteration approach}.
\bjtitle{IEEE Transactions on Neural Networks and Learning Systems}
\bvolume{31}(\bissue{2}),
\bfpage{549}--\blpage{558}
(\byear{2020})
\end{barticle}
\endbibitem

\bibitem[\protect\citeauthoryear{Luo et~al.}{2014}]{DataHJB}
\begin{barticle}
\bauthor{\bsnm{Luo}, \binits{B.}},
\bauthor{\bsnm{Wu}, \binits{H.-N.}},
\bauthor{\bsnm{Huang}, \binits{T.}},
\bauthor{\bsnm{Liu}, \binits{D.}}:
\batitle{Data-based approximate policy iteration for affine nonlinear
  continuous-time optimal control design}.
\bjtitle{Automatica}
\bvolume{50}(\bissue{12}),
\bfpage{3281}--\blpage{3290}
(\byear{2014})
\end{barticle}
\endbibitem

\bibitem[\protect\citeauthoryear{Pakkhesal and Shamaghdari}{}]{SOS}
\begin{botherref}
\oauthor{\bsnm{Pakkhesal}, \binits{S.}},
\oauthor{\bsnm{Shamaghdari}, \binits{S.}}:
Sum-of-squares-based policy iteration for suboptimal control of polynomial
  time-varying systems.
Asian Journal of Control
\textbf{n/a}
\end{botherref}
\endbibitem

\bibitem[\protect\citeauthoryear{Yazdani et~al.}{2020}]{safetySOS}
\begin{barticle}
\bauthor{\bsnm{Yazdani}, \binits{N.}},
\bauthor{\bsnm{Moghaddam}, \binits{R.}},
\bauthor{\bsnm{Kiumarsi}, \binits{B.}},
\bauthor{\bsnm{Modares}, \binits{H.}}:
\batitle{A safety-certified policy iteration algorithm for control of
  constrained nonlinear systems}.
\bjtitle{IEEE Control Systems Letters}
\bvolume{4}(\bissue{3}),
\bfpage{686}--\blpage{691}
(\byear{2020})
\end{barticle}
\endbibitem

\bibitem[\protect\citeauthoryear{Tonon et~al.}{2017}]{SemiLagranigian}
\begin{bbook}
\bauthor{\bsnm{Tonon}, \binits{D.}},
\bauthor{\bsnm{Aronna}, \binits{M.}},
\bauthor{\bsnm{Kalise}, \binits{D.}}:
\bbtitle{Optimal Control: Novel Directions and Applications}.
\bpublisher{Springer},
\blocation{International}
(\byear{2017})
\end{bbook}
\endbibitem

\bibitem[\protect\citeauthoryear{Debrabant and
  Jakobsen}{2014}]{SemiLagrangianStochastic}
\begin{bchapter}
\bauthor{\bsnm{Debrabant}, \binits{K.}},
\bauthor{\bsnm{Jakobsen}, \binits{E.}}:
\bctitle{Semi-{L}agrangian schemes for linear and fully non-linear
  {H}amilton-{J}acobi-{B}ellman equations}.
In: \bbtitle{Hyperbolic Problems: Theory, Numerics, Applications},
pp. \bfpage{483}--\blpage{490}.
\bpublisher{Springer},
\blocation{International}
(\byear{2014})
\end{bchapter}
\endbibitem

\bibitem[\protect\citeauthoryear{Falcone}{1987}]{Falcone1987}
\begin{barticle}
\bauthor{\bsnm{Falcone}, \binits{M.}}:
\batitle{A numerical approach to the infinite horizon problem of deterministic
  control theory}.
\bjtitle{Applied Mathematics and Optimization}
\bvolume{15}(\bissue{1}),
\bfpage{1}--\blpage{13}
(\byear{1987})
\end{barticle}
\endbibitem

\bibitem[\protect\citeauthoryear{Falcone et~al.}{1994}]{FALCONESplitting}
\begin{barticle}
\bauthor{\bsnm{Falcone}, \binits{M.}},
\bauthor{\bsnm{Lanucara}, \binits{P.}},
\bauthor{\bsnm{Seghini}, \binits{A.}}:
\batitle{A splitting algorithm for {H}amilton-{J}acobi-{B}ellman equations}.
\bjtitle{Applied Numerical Mathematics}
\bvolume{15}(\bissue{2}),
\bfpage{207}--\blpage{218}
(\byear{1994})
\end{barticle}
\endbibitem

\bibitem[\protect\citeauthoryear{Kafash et~al.}{2013}]{VIM}
\begin{barticle}
\bauthor{\bsnm{Kafash}, \binits{B.}},
\bauthor{\bsnm{Delavarkhalafi}, \binits{A.}},
\bauthor{\bsnm{Karbassi}, \binits{S.M.}}:
\batitle{Application of variational iteration method for
  {H}amilton-{J}acobi-{B}ellman}.
\bjtitle{Applied Mathematical Modelling}
\bvolume{37}(\bissue{6}),
\bfpage{3917}--\blpage{3928}
(\byear{2013})
\end{barticle}
\endbibitem

\bibitem[\protect\citeauthoryear{Alla and Saluzzi}{2020}]{ALLA2020192}
\begin{barticle}
\bauthor{\bsnm{Alla}, \binits{A.}},
\bauthor{\bsnm{Saluzzi}, \binits{L.}}:
\batitle{{A HJB-POD approach for the control of nonlinear PDEs on a tree
  structure}}.
\bjtitle{Applied Numerical Mathematics}
\bvolume{155},
\bfpage{192}--\blpage{207}
(\byear{2020})
\end{barticle}
\endbibitem

\bibitem[\protect\citeauthoryear{Akian et~al.}{2009}]{maxplus_det}
\begin{bchapter}
\bauthor{\bsnm{Akian}, \binits{M.}},
\bauthor{\bsnm{Gaubert}, \binits{S.}},
\bauthor{\bsnm{Lakhoua}, \binits{A.}}:
\bctitle{Convergence analysis of the max-plus finite element method for solving
  deterministic optimal control problems}.
In: \bbtitle{Proceedings of the IEEE Conference on Decision and Control},
pp. \bfpage{927}--\blpage{934}.
\bpublisher{IEEE},
\blocation{NY}
(\byear{2009})
\end{bchapter}
\endbibitem

\bibitem[\protect\citeauthoryear{Akian and Fodjo}{2018}]{maxplus_stoch}
\begin{bbook}
\bauthor{\bsnm{Akian}, \binits{M.}},
\bauthor{\bsnm{Fodjo}, \binits{E.}}:
\bbtitle{Probabilistic Max-Plus Schemes for Solving Hamilton-Jacobi-Bellman
  Equations},
pp. \bfpage{183}--\blpage{209}.
\bpublisher{Springer},
\blocation{International}
(\byear{2018})
\end{bbook}
\endbibitem

\bibitem[\protect\citeauthoryear{Pontryagin et~al.}{1962}]{Pontryagin}
\begin{bbook}
\bauthor{\bsnm{Pontryagin}, \binits{L.}},
\bauthor{\bsnm{Boltyanskii}, \binits{V.}},
\bauthor{\bsnm{Gamkrelidze}, \binits{R.}},
\bauthor{\bsnm{Mishchenko}, \binits{E.}}:
\bbtitle{{The Mathematical Theory of Optimal Processes}}.
\bpublisher{Translated from the Russian by K. N. Trirogoff; edited by L.
  W.Neustadt. Wiley},
\blocation{New York, NY}
(\byear{1962})
\end{bbook}
\endbibitem

\bibitem[\protect\citeauthoryear{Beeler et~al.}{2000}]{beeler2000feedback}
\begin{barticle}
\bauthor{\bsnm{Beeler}, \binits{S.}},
\bauthor{\bsnm{Tran}, \binits{H.}},
\bauthor{\bsnm{Banks}, \binits{H.}}:
\batitle{Feedback control methodologies for nonlinear systems}.
\bjtitle{Journal of optimization theory and applications}
\bvolume{107}(\bissue{1}),
\bfpage{1}--\blpage{33}
(\byear{2000})
\end{barticle}
\endbibitem

\bibitem[\protect\citeauthoryear{Kang and Wilcox}{2017}]{kang2017mitigating}
\begin{barticle}
\bauthor{\bsnm{Kang}, \binits{W.}},
\bauthor{\bsnm{Wilcox}, \binits{L.}}:
\batitle{Mitigating the curse of dimensionality: sparse grid characteristics
  method for optimal feedback control and hjb equations}.
\bjtitle{Computational Optimization and Applications}
\bvolume{68}(\bissue{2}),
\bfpage{289}--\blpage{315}
(\byear{2017})
\end{barticle}
\endbibitem

\bibitem[\protect\citeauthoryear{Nakamura-Zimmerer
  et~al.}{2021}]{nakamura2019adaptive}
\begin{barticle}
\bauthor{\bsnm{Nakamura-Zimmerer}, \binits{T.}},
\bauthor{\bsnm{Gong}, \binits{Q.}},
\bauthor{\bsnm{Kang}, \binits{W.}}:
\batitle{Adaptive deep learning for high-dimensional hamilton--jacobi--bellman
  equations}.
\bjtitle{SIAM Journal on Scientific Computing}
\bvolume{43}(\bissue{2}),
\bfpage{1221}--\blpage{1247}
(\byear{2021})
\end{barticle}
\endbibitem

\bibitem[\protect\citeauthoryear{Azmi et~al.}{2021}]{azmi2020optimal}
\begin{barticle}
\bauthor{\bsnm{Azmi}, \binits{B.}},
\bauthor{\bsnm{Kalise}, \binits{D.}},
\bauthor{\bsnm{Kunisch}, \binits{K.}}:
\batitle{Optimal feedback law recovery by gradient-augmented sparse polynomial
  regression}.
\bjtitle{Journal of Machine Learning Research}
\bvolume{22},
\bfpage{1}--\blpage{32}
(\byear{2021})
\end{barticle}
\endbibitem

\bibitem[\protect\citeauthoryear{Vapnik}{1992}]{vapnik1992principles}
\begin{bchapter}
\bauthor{\bsnm{Vapnik}, \binits{V.}}:
\bctitle{Principles of risk minimization for learning theory}.
In: \bbtitle{Advances in Neural Information Processing Systems},
pp. \bfpage{831}--\blpage{838}
(\byear{1992})
\end{bchapter}
\endbibitem

\bibitem[\protect\citeauthoryear{Steinwart and
  Christmann}{2008}]{steinwart2008support}
\begin{bbook}
\bauthor{\bsnm{Steinwart}, \binits{I.}},
\bauthor{\bsnm{Christmann}, \binits{A.}}:
\bbtitle{Support Vector Machines}.
\bpublisher{Springer},
\blocation{Berlin}
(\byear{2008})
\end{bbook}
\endbibitem

\bibitem[\protect\citeauthoryear{Hackbusch}{2012}]{hackbusch-2012}
\begin{bbook}
\bauthor{\bsnm{Hackbusch}, \binits{W.}}:
\bbtitle{Tensor Spaces And Numerical Tensor Calculus}.
\bpublisher{Springer},
\blocation{Berlin}
(\byear{2012})
\end{bbook}
\endbibitem

\bibitem[\protect\citeauthoryear{Oseledets and Tyrtyshnikov}{2009}]{Oseledets}
\begin{barticle}
\bauthor{\bsnm{Oseledets}, \binits{I.}},
\bauthor{\bsnm{Tyrtyshnikov}, \binits{E.}}:
\batitle{Breaking the curse of dimensionality, or how to use {SVD} in many
  dimensions}.
\bjtitle{SIAM J. Sci. Comput.}
\bvolume{31},
\bfpage{3744}--\blpage{3759}
(\byear{2009})
\end{barticle}
\endbibitem

\bibitem[\protect\citeauthoryear{Oseledets}{2011}]{osel-tt-2011}
\begin{barticle}
\bauthor{\bsnm{Oseledets}, \binits{I.V.}}:
\batitle{Tensor-train decomposition}.
\bjtitle{SIAM J. Sci. Comput.}
\bvolume{33}(\bissue{5}),
\bfpage{2295}--\blpage{2317}
(\byear{2011})
\end{barticle}
\endbibitem

\bibitem[\protect\citeauthoryear{Khoromskij}{2011}]{Khoromskij-book}
\begin{barticle}
\bauthor{\bsnm{Khoromskij}, \binits{B.N.}}:
\batitle{{Tensors-structured numerical methods in scientific computing : survey
  on recent advances}}.
\bjtitle{Chemometrics and intelligent laboratory systems}
\bvolume{110}(\bissue{1}),
\bfpage{1}--\blpage{19}
(\byear{2011})
\end{barticle}
\endbibitem

\bibitem[\protect\citeauthoryear{Hackbusch and Schneider}{2014}]{Hackbusch2014}
\begin{bbook}
\bauthor{\bsnm{Hackbusch}, \binits{W.}},
\bauthor{\bsnm{Schneider}, \binits{R.}}:
\bbtitle{Tensor Spaces and Hierarchical Tensor Representations},
pp. \bfpage{237}--\blpage{261}.
\bpublisher{Springer},
\blocation{Cham}
(\byear{2014})
\end{bbook}
\endbibitem

\bibitem[\protect\citeauthoryear{Bachmayr
  et~al.}{2016}]{Bachmayr-Uschmajew-Schneider}
\begin{barticle}
\bauthor{\bsnm{Bachmayr}, \binits{M.}},
\bauthor{\bsnm{Schneider}, \binits{R.}},
\bauthor{\bsnm{Uschmajew}, \binits{A.}}:
\batitle{Tensor networks and hierarchical tensors for the solution of
  high-dimensional partial differential equations}.
\bjtitle{Found. Comput. Math.}
\bvolume{16}(\bissue{6}),
\bfpage{1423}--\blpage{1472}
(\byear{2016})
\end{barticle}
\endbibitem

\bibitem[\protect\citeauthoryear{Szalay et~al.}{2015}]{Legeza-Schneider}
\begin{barticle}
\bauthor{\bsnm{Szalay}, \binits{S.}},
\bauthor{\bsnm{Pfeffer}, \binits{M.}},
\bauthor{\bsnm{Murg}, \binits{V.}},
\bauthor{\bsnm{Barcza}, \binits{G.}},
\bauthor{\bsnm{Verstraete}, \binits{F.}},
\bauthor{\bsnm{Schneider}, \binits{R.}},
\bauthor{\bsnm{Legeza}}:
\batitle{{Tensor product methods and entanglement optimization for ab initio
  quantum chemistry}}.
\bjtitle{International j. of quantum chemistry}
\bvolume{115}(\bissue{19}),
\bfpage{1342}--\blpage{1391}
(\byear{2015})
\end{barticle}
\endbibitem

\bibitem[\protect\citeauthoryear{Hackbusch}{2014}]{Hackbusch-Acta}
\begin{barticle}
\bauthor{\bsnm{Hackbusch}, \binits{W.}}:
\batitle{{Numerical tensor calculus}}.
\bjtitle{Acta numerica}
\bvolume{23},
\bfpage{651}--\blpage{742}
(\byear{2014})
\end{barticle}
\endbibitem

\bibitem[\protect\citeauthoryear{Bachmayr et~al.}{2017}]{Dahmen-Bachmayr}
\begin{barticle}
\bauthor{\bsnm{Bachmayr}, \binits{M.}},
\bauthor{\bsnm{Cohen}, \binits{A.}},
\bauthor{\bsnm{Dahmen}, \binits{W.}}:
\batitle{{Parametric PDEs: sparse or low-rank approximations?}}
\bjtitle{IMA J. of Numerical Analysis}
\bvolume{38}(\bissue{4}),
\bfpage{1661}--\blpage{1708}
(\byear{2017})
\end{barticle}
\endbibitem

\bibitem[\protect\citeauthoryear{Dolgov et~al.}{2021}]{DKK21}
\begin{barticle}
\bauthor{\bsnm{Dolgov}, \binits{S.}},
\bauthor{\bsnm{Kalise}, \binits{D.}},
\bauthor{\bsnm{Kunisch}, \binits{K.K.}}:
\batitle{Tensor {Decomposition} {Methods} for {High}-dimensional
  {Hamilton}--{Jacobi}--{Bellman} {Equations}}.
\bjtitle{{SIAM} Journal on Scientific Computing}
\bvolume{43}(\bissue{3}),
\bfpage{1625}--\blpage{1650}
(\byear{2021})
\end{barticle}
\endbibitem

\bibitem[\protect\citeauthoryear{Oster et~al.}{2022}]{oster2}
\begin{barticle}
\bauthor{\bsnm{Oster}, \binits{M.}},
\bauthor{\bsnm{Sallandt}, \binits{L.}},
\bauthor{\bsnm{Schneider}, \binits{R.}}:
\batitle{Approximating optimal feedback controllers of finite horizon control
  problems using hierarchical tensor formats}.
\bjtitle{SIAM Journal on Scientific Computing}
\bvolume{44}(\bissue{3}),
\bfpage{746}--\blpage{770}
(\byear{2022})
\end{barticle}
\endbibitem

\bibitem[\protect\citeauthoryear{Oster et~al.}{2019}]{oster3}
\begin{botherref}
\oauthor{\bsnm{Oster}, \binits{M.}},
\oauthor{\bsnm{Sallandt}, \binits{L.}},
\oauthor{\bsnm{Schneider}, \binits{R.}}:
Approximating the stationary bellman equation by hierarchical tensor products.
Journal of Computational Mathematics
(2019)
\end{botherref}
\endbibitem

\bibitem[\protect\citeauthoryear{Stefansson and Leong}{2016}]{stefansson}
\begin{bchapter}
\bauthor{\bsnm{Stefansson}, \binits{E.}},
\bauthor{\bsnm{Leong}, \binits{Y.}}:
\bctitle{Sequential alternating least squares for solving high dimensional
  linear hamilton-jacobi-bellman equation}.
In: \bbtitle{2016 IEEE/RSJ International Conference on Intelligent Robots and
  Systems (IROS)},
pp. \bfpage{3757}--\blpage{3764}
(\byear{2016})
\end{bchapter}
\endbibitem

\bibitem[\protect\citeauthoryear{Horowitz and
  Burdick}{2014}]{horowitz2014linear}
\begin{bchapter}
\bauthor{\bsnm{Horowitz}, \binits{A.} \bsuffix{M.and~Damle}},
\bauthor{\bsnm{Burdick}, \binits{J.}}:
\bctitle{Linear {Hamilton Jacobi Bellman} equations in high dimensions}.
In: \bbtitle{53rd IEEE Conference on Decision and Control},
pp. \bfpage{5880}--\blpage{5887}
(\byear{2014}).
\bcomment{IEEE}
\end{bchapter}
\endbibitem

\bibitem[\protect\citeauthoryear{Fackeldey et~al.}{2022}]{oster1}
\begin{barticle}
\bauthor{\bsnm{Fackeldey}, \binits{K.}},
\bauthor{\bsnm{Oster}, \binits{M.}},
\bauthor{\bsnm{Sallandt}, \binits{L.}},
\bauthor{\bsnm{Schneider}, \binits{R.}}:
\batitle{Approximative policy iteration for exit time feedback control problems
  driven by stochastic differential equations using tensor train format}.
\bjtitle{Multiscale Modeling \& Simulation}
\bvolume{20}(\bissue{1}),
\bfpage{379}--\blpage{403}
(\byear{2022})
\end{barticle}
\endbibitem

\bibitem[\protect\citeauthoryear{Gorodetsky et~al.}{2018}]{gorodetsky2018high}
\begin{barticle}
\bauthor{\bsnm{Gorodetsky}, \binits{A.}},
\bauthor{\bsnm{Karaman}, \binits{S.}},
\bauthor{\bsnm{Marzouk}, \binits{Y.}}:
\batitle{High-dimensional stochastic optimal control using continuous tensor
  decompositions}.
\bjtitle{The International Journal of Robotics Research}
\bvolume{37}(\bissue{2-3}),
\bfpage{340}--\blpage{377}
(\byear{2018})
\end{barticle}
\endbibitem

\bibitem[\protect\citeauthoryear{Dolgov et~al.}{2023}]{dolgov2022data}
\begin{barticle}
\bauthor{\bsnm{Dolgov}, \binits{S.}},
\bauthor{\bsnm{Kalise}, \binits{D.}},
\bauthor{\bsnm{Saluzzi}, \binits{L.}}:
\batitle{Data-driven tensor train gradient cross approximation for
  hamilton--jacobi--bellman equations}.
\bjtitle{SIAM Journal on Scientific Computing}
\bvolume{45}(\bissue{5}),
\bfpage{2153}--\blpage{2184}
(\byear{2023})
\end{barticle}
\endbibitem

\bibitem[\protect\citeauthoryear{Götte et~al.}{2021}]{Micha}
\begin{botherref}
\oauthor{\bsnm{Götte}, \binits{M.}},
\oauthor{\bsnm{Schneider}, \binits{R.}},
\oauthor{\bsnm{Trunschke}, \binits{P.}}:
A block-sparse tensor train format for sample-efficient high-dimensional
  polynomial regression.
Frontiers in Applied Mathematics and Statistics
\textbf{7}
(2021)
\end{botherref}
\endbibitem

\bibitem[\protect\citeauthoryear{Oseledets and
  Tyrtyshnikov}{2010}]{ot-ttcross-2010}
\begin{barticle}
\bauthor{\bsnm{Oseledets}, \binits{I.V.}},
\bauthor{\bsnm{Tyrtyshnikov}, \binits{E.E.}}:
\batitle{{TT-cross} approximation for multidimensional arrays}.
\bjtitle{Linear Algebra Appl.}
\bvolume{432}(\bissue{1}),
\bfpage{70}--\blpage{88}
(\byear{2010})
\end{barticle}
\endbibitem

\bibitem[\protect\citeauthoryear{Savostyanov and
  Oseledets}{2011}]{so-dmrgi-2011proc}
\begin{bchapter}
\bauthor{\bsnm{Savostyanov}, \binits{D.V.}},
\bauthor{\bsnm{Oseledets}, \binits{I.V.}}:
\bctitle{Fast adaptive interpolation of multi-dimensional arrays in tensor
  train format}.
In: \bbtitle{Proceedings of 7th International Workshop on Multidimensional
  Systems (nDS)}.
\bpublisher{IEEE},
\blocation{NY}
(\byear{2011})
\end{bchapter}
\endbibitem

\bibitem[\protect\citeauthoryear{Grasedyck
  et~al.}{2015}]{grasedyck-par-cross-2015}
\begin{barticle}
\bauthor{\bsnm{Grasedyck}, \binits{L.}},
\bauthor{\bsnm{Kriemann}, \binits{R.}},
\bauthor{\bsnm{L{\"o}bbert}, \binits{C.}},
\bauthor{\bsnm{N{\"a}gel}, \binits{A.}},
\bauthor{\bsnm{Wittum}, \binits{G.}},
\bauthor{\bsnm{Xylouris}, \binits{K.}}:
\batitle{Parallel tensor sampling in the hierarchical {Tucker} format}.
\bjtitle{Computing and Visualization in Science}
\bvolume{17}(\bissue{2}),
\bfpage{67}--\blpage{78}
(\byear{2015})
\end{barticle}
\endbibitem

\bibitem[\protect\citeauthoryear{Savostyanov}{2014}]{sav-qott-2014}
\begin{barticle}
\bauthor{\bsnm{Savostyanov}, \binits{D.V.}}:
\batitle{Quasioptimality of maximum--volume cross interpolation of tensors}.
\bjtitle{Linear Algebra Appl.}
\bvolume{458},
\bfpage{217}--\blpage{244}
(\byear{2014})
\end{barticle}
\endbibitem

\bibitem[\protect\citeauthoryear{Wendland}{2005}]{wendland2005scattered}
\begin{bbook}
\bauthor{\bsnm{Wendland}, \binits{H.}}:
\bbtitle{Scattered {D}ata {A}pproximation}.
\bsertitle{Cambridge Monographs on Applied and Computational Mathematics},
vol. \bseriesno{17}.
\bpublisher{Cambridge University Press},
\blocation{Cambridge}
(\byear{2005})
\end{bbook}
\endbibitem

\bibitem[\protect\citeauthoryear{Berner et~al.}{2022}]{DeepLearningGitta}
\begin{bchapter}
\bauthor{\bsnm{Berner}, \binits{J.}},
\bauthor{\bsnm{Grohs}, \binits{P.}},
\bauthor{\bsnm{Kutyniok}, \binits{G.}},
\bauthor{\bsnm{Petersen}, \binits{P.}}:
\bctitle{The modern mathematics of deep learning}.
In: \bbtitle{Mathematical Aspects of Deep Learning},
pp. \bfpage{1}--\blpage{111}.
\bpublisher{Cambridge University Press},
\blocation{Cambridge}
(\byear{2022})
\end{bchapter}
\endbibitem

\bibitem[\protect\citeauthoryear{DeVore et~al.}{2021}]{DHP}
\begin{barticle}
\bauthor{\bsnm{DeVore}, \binits{R.A.}},
\bauthor{\bsnm{Hanin}, \binits{B.}},
\bauthor{\bsnm{Petrova}, \binits{G.}}:
\batitle{Neural network approximation}.
\bjtitle{Acta Numerica}
\bvolume{30},
\bfpage{327}--\blpage{444}
(\byear{2021})
\end{barticle}
\endbibitem

\bibitem[\protect\citeauthoryear{E et~al.}{2020}]{DeepLearningWeinan}
\begin{botherref}
\oauthor{\bsnm{E}, \binits{W.}},
\oauthor{\bsnm{Ma}, \binits{C.}},
\oauthor{\bsnm{Wojtowytsch}, \binits{S.}},
\oauthor{\bsnm{Wu}, \binits{L.}}:
Towards a Mathematical Understanding of Neural Network-Based Machine Learning:
  what we know and what we don't.
arXiv
(2020)
\end{botherref}
\endbibitem

\bibitem[\protect\citeauthoryear{Higham and Higham}{2019}]{DeepLearningHigham}
\begin{barticle}
\bauthor{\bsnm{Higham}, \binits{C.F.}},
\bauthor{\bsnm{Higham}, \binits{D.J.}}:
\batitle{Deep learning: An introduction for applied mathematicians}.
\bjtitle{SIAM Review}
\bvolume{61}(\bissue{4}),
\bfpage{860}--\blpage{891}
(\byear{2019})
\end{barticle}
\endbibitem

\bibitem[\protect\citeauthoryear{Pak and Kim}{2017}]{ImageRec}
\begin{bchapter}
\bauthor{\bsnm{Pak}, \binits{M.}},
\bauthor{\bsnm{Kim}, \binits{S.}}:
\bctitle{A review of deep learning in image recognition}.
In: \bbtitle{2017 4th International Conference on Computer Applications and
  Information Processing Technology (CAIPT)},
pp. \bfpage{1}--\blpage{3}
(\byear{2017})
\end{bchapter}
\endbibitem

\bibitem[\protect\citeauthoryear{Beck et~al.}{2023}]{DeepLearningPDE}
\begin{barticle}
\bauthor{\bsnm{Beck}, \binits{C.}},
\bauthor{\bsnm{Hutzenthaler}, \binits{M.}},
\bauthor{\bsnm{Jentzen}, \binits{A.}},
\bauthor{\bsnm{Kuckuck}, \binits{B.}}:
\batitle{An overview on deep learning-based approximation methods for partial
  differential equations}.
\bjtitle{Discrete and Continuous Dynamical Systems - B}
\bvolume{28}(\bissue{6}),
\bfpage{3697}--\blpage{3746}
(\byear{2023})
\end{barticle}
\endbibitem

\bibitem[\protect\citeauthoryear{{Kunisch, Karl} and {Walter,
  Daniel}}{2021}]{Kunisch}
\begin{barticle}
\bauthor{\bsnm{{Kunisch, Karl}}},
\bauthor{\bsnm{{Walter, Daniel}}}:
\batitle{Semiglobal optimal feedback stabilization of autonomous systems via
  deep neural network approximation}.
\bjtitle{ESAIM: COCV}
\bvolume{27},
\bfpage{16}
(\byear{2021})
\end{barticle}
\endbibitem

\bibitem[\protect\citeauthoryear{Kunisch and Walter}{2023}]{Kunisch2}
\begin{botherref}
\oauthor{\bsnm{Kunisch}, \binits{K.}},
\oauthor{\bsnm{Walter}, \binits{D.}}:
Optimal feedback control of dynamical systems via value-function approximation.
arXiv preprint arXiv:2302.13122
(2023)
\end{botherref}
\endbibitem

\bibitem[\protect\citeauthoryear{Darbon et~al.}{2020}]{darbon2020overcoming}
\begin{barticle}
\bauthor{\bsnm{Darbon}, \binits{J.}},
\bauthor{\bsnm{Langlois}, \binits{G.P.}},
\bauthor{\bsnm{Meng}, \binits{T.}}:
\batitle{Overcoming the curse of dimensionality for some hamilton--jacobi
  partial differential equations via neural network architectures}.
\bjtitle{Research in the Mathematical Sciences}
\bvolume{7}(\bissue{3}),
\bfpage{1}--\blpage{50}
(\byear{2020})
\end{barticle}
\endbibitem

\bibitem[\protect\citeauthoryear{N{\"u}sken and
  Richter}{2021}]{nusken2020solving}
\begin{barticle}
\bauthor{\bsnm{N{\"u}sken}, \binits{N.}},
\bauthor{\bsnm{Richter}, \binits{L.}}:
\batitle{Solving high-dimensional hamilton--jacobi--bellman pdes using neural
  networks: perspectives from the theory of controlled diffusions and measures
  on path space}.
\bjtitle{Partial Differential Equations and Applications}
\bvolume{2}(\bissue{4}),
\bfpage{1}--\blpage{48}
(\byear{2021})
\end{barticle}
\endbibitem

\bibitem[\protect\citeauthoryear{Ito et~al.}{2020}]{ito2020neural}
\begin{botherref}
\oauthor{\bsnm{Ito}, \binits{K.}},
\oauthor{\bsnm{Reisinger}, \binits{C.}},
\oauthor{\bsnm{Zhang}, \binits{Y.}}:
A neural network-based policy iteration algorithm with global $h^{2}$
  -superlinear convergence for stochastic games on domains.
Foundations of Computational Mathematics,
1--44
(2020)
\end{botherref}
\endbibitem

\bibitem[\protect\citeauthoryear{Demo et~al.}{2023}]{DEMO2023383}
\begin{barticle}
\bauthor{\bsnm{Demo}, \binits{N.}},
\bauthor{\bsnm{Strazzullo}, \binits{M.}},
\bauthor{\bsnm{Rozza}, \binits{G.}}:
\batitle{An extended physics informed neural network for preliminary analysis
  of parametric optimal control problems}.
\bjtitle{Computers \& Mathematics with Applications}
\bvolume{143},
\bfpage{383}--\blpage{396}
(\byear{2023})
\end{barticle}
\endbibitem

\bibitem[\protect\citeauthoryear{Han et~al.}{2018}]{Han_Jentzen_E_2018}
\begin{barticle}
\bauthor{\bsnm{Han}, \binits{J.}},
\bauthor{\bsnm{Jentzen}, \binits{A.}},
\bauthor{\bsnm{E}, \binits{W.}}:
\batitle{Solving high-dimensional partial differential equations using deep
  learning}.
\bjtitle{Proceedings of the National Academy of Sciences}
\bvolume{115}(\bissue{34}),
\bfpage{8505}--\blpage{8510}
(\byear{2018})
\end{barticle}
\endbibitem

\bibitem[\protect\citeauthoryear{Meng et~al.}{2022}]{sympocnet}
\begin{botherref}
\oauthor{\bsnm{Meng}, \binits{T.}},
\oauthor{\bsnm{Zhang}, \binits{Z.}},
\oauthor{\bsnm{Darbon}, \binits{J.}},
\oauthor{\bsnm{Karniadakis}, \binits{G.E.}}:
SympOCnet: Solving optimal control problems with applications to
  high-dimensional multi-agent path planning problems.
arXiv.
OPTdoi: 10.48550/ARXIV.2201.05475
(2022)
\end{botherref}
\endbibitem

\bibitem[\protect\citeauthoryear{Zhou et~al.}{2021}]{Zhou_2021}
\begin{barticle}
\bauthor{\bsnm{Zhou}, \binits{M.}},
\bauthor{\bsnm{Han}, \binits{J.}},
\bauthor{\bsnm{Lu}, \binits{J.}}:
\batitle{Actor-critic method for high dimensional static
  hamilton--jacobi--bellman partial differential equations based on neural
  networks}.
\bjtitle{{SIAM} Journal on Scientific Computing}
\bvolume{43}(\bissue{6}),
\bfpage{4043}--\blpage{4066}
(\byear{2021})
\end{barticle}
\endbibitem

\bibitem[\protect\citeauthoryear{Onken et~al.}{2021}]{Onken2021}
\begin{bchapter}
\bauthor{\bsnm{Onken}, \binits{D.}},
\bauthor{\bsnm{Nurbekyan}, \binits{L.}},
\bauthor{\bsnm{Li}, \binits{X.}},
\bauthor{\bsnm{Fung}, \binits{S.W.}},
\bauthor{\bsnm{Osher}, \binits{S.}},
\bauthor{\bsnm{Ruthotto}, \binits{L.}}:
\bctitle{A neural network approach applied to multi-agent optimal control}.
In: \bbtitle{2021 European Control Conference ({ECC})}.
\bpublisher{{IEEE}},
\blocation{NY}
(\byear{2021})
\end{bchapter}
\endbibitem

\bibitem[\protect\citeauthoryear{Ruthotto et~al.}{2020}]{ruthotto2020machine}
\begin{barticle}
\bauthor{\bsnm{Ruthotto}, \binits{L.}},
\bauthor{\bsnm{Osher}, \binits{S.J.}},
\bauthor{\bsnm{Li}, \binits{W.}},
\bauthor{\bsnm{Nurbekyan}, \binits{L.}},
\bauthor{\bsnm{Fung}, \binits{S.W.}}:
\batitle{A machine learning framework for solving high-dimensional mean field
  game and mean field control problems}.
\bjtitle{Proceedings of the National Academy of Sciences}
\bvolume{117}(\bissue{17}),
\bfpage{9183}--\blpage{9193}
(\byear{2020})
\end{barticle}
\endbibitem

\bibitem[\protect\citeauthoryear{Albi et~al.}{2022}]{ABK21}
\begin{barticle}
\bauthor{\bsnm{Albi}, \binits{G.}},
\bauthor{\bsnm{Bicego}, \binits{S.}},
\bauthor{\bsnm{Kalise}, \binits{D.}}:
\batitle{Gradient-augmented {Supervised} {Learning} of {Optimal} {Feedback}
  {Laws} {Using} {State}-{Dependent} {Riccat} {Equations}}.
\bjtitle{IEEE Control Systems Letters}
\bvolume{6},
\bfpage{836}--\blpage{841}
(\byear{2022})
\end{barticle}
\endbibitem

\bibitem[\protect\citeauthoryear{Gr{\"u}ne}{2020}]{grune2020computing}
\begin{botherref}
\oauthor{\bsnm{Gr{\"u}ne}, \binits{L.}}:
Computing lyapunov functions using deep neural networks.
arXiv preprint arXiv:2005.08965
(2020)
\end{botherref}
\endbibitem

\bibitem[\protect\citeauthoryear{Kunisch et~al.}{2021}]{Kunisch3}
\begin{barticle}
\bauthor{\bsnm{Kunisch}, \binits{K.}},
\bauthor{\bsnm{Rodrigues}, \binits{S.S.}},
\bauthor{\bsnm{Walter}, \binits{D.}}:
\batitle{Learning an optimal feedback operator semiglobally stabilizing
  semilinear parabolic equations}.
\bjtitle{Applied Mathematics {\&} Optimization}
\bvolume{84}(\bissue{1}),
\bfpage{277}--\blpage{318}
(\byear{2021})
\end{barticle}
\endbibitem

\bibitem[\protect\citeauthoryear{Kunisch et~al.}{2022}]{Kunisch4}
\begin{botherref}
\oauthor{\bsnm{Kunisch}, \binits{K.}},
\oauthor{\bsnm{Vásquez-Varas}, \binits{D.}},
\oauthor{\bsnm{Walter}, \binits{D.}}:
Learning Optimal Feedback Operators and their Polynomial Approximation.
arXiv
(2022)
\end{botherref}
\endbibitem

\bibitem[\protect\citeauthoryear{Azmi et~al.}{2021}]{azmi2021optimal}
\begin{barticle}
\bauthor{\bsnm{Azmi}, \binits{B.}},
\bauthor{\bsnm{Kalise}, \binits{D.}},
\bauthor{\bsnm{Kunisch}, \binits{K.}}:
\batitle{Optimal feedback law recovery by gradient-augmented sparse polynomial
  regression}.
\bjtitle{Journal of Machine Learning Research}
\bvolume{22}(\bissue{48}),
\bfpage{1}--\blpage{32}
(\byear{2021})
\end{barticle}
\endbibitem

\bibitem[\protect\citeauthoryear{{\c{C}}imen}{2008}]{ccimen2008state}
\begin{barticle}
\bauthor{\bsnm{{\c{C}}imen}, \binits{T.}}:
\batitle{State-dependent {Riccati equation (SDRE)} control: a survey}.
\bjtitle{IFAC Proceedings Volumes}
\bvolume{41}(\bissue{2}),
\bfpage{3761}--\blpage{3775}
(\byear{2008})
\end{barticle}
\endbibitem

\bibitem[\protect\citeauthoryear{Alla et~al.}{2021}]{allasdre}
\begin{botherref}
\oauthor{\bsnm{Alla}, \binits{A.}},
\oauthor{\bsnm{Kalise}, \binits{D.}},
\oauthor{\bsnm{Simoncini}, \binits{V.}}:
State-dependent {Riccati} equation feedback stabilization for nonlinear {PDEs}.
OPTdoi: 10.48550/ARXIV.2106.07163
(2021)
\end{botherref}
\endbibitem

\bibitem[\protect\citeauthoryear{Banks et~al.}{2007}]{Banks_Lewis_Tran_2007}
\begin{barticle}
\bauthor{\bsnm{Banks}, \binits{H.T.}},
\bauthor{\bsnm{Lewis}, \binits{B.M.}},
\bauthor{\bsnm{Tran}, \binits{H.T.}}:
\batitle{Nonlinear feedback controllers and compensators: a state-dependent
  {Riccati} equation approach}.
\bjtitle{Computational Optimization and Applications}
\bvolume{37}(\bissue{2}),
\bfpage{177}--\blpage{218}
(\byear{2007})
\end{barticle}
\endbibitem

\bibitem[\protect\citeauthoryear{Rohrbach et~al.}{2022}]{rdgs-tt-gauss-2022}
\begin{barticle}
\bauthor{\bsnm{Rohrbach}, \binits{P.B.}},
\bauthor{\bsnm{Dolgov}, \binits{S.}},
\bauthor{\bsnm{Grasedyck}, \binits{L.}},
\bauthor{\bsnm{Scheichl}, \binits{R.}}:
\batitle{Rank bounds for approximating {Gaussian} densities in the
  {Tensor-Train} format}.
\bjtitle{SIAM/ASA Journal on Uncertainty Quantification}
\bvolume{10}(\bissue{3}),
\bfpage{1191}--\blpage{1224}
(\byear{2022})
\end{barticle}
\endbibitem

\bibitem[\protect\citeauthoryear{Holtz et~al.}{2012}]{ALS}
\begin{barticle}
\bauthor{\bsnm{Holtz}, \binits{S.}},
\bauthor{\bsnm{Rohwedder}, \binits{T.}},
\bauthor{\bsnm{Schneider}, \binits{R.}}:
\batitle{The alternating linear scheme for tensor optimization in the tensor
  train format}.
\bjtitle{SIAM J. Sci. Comput.}
\bvolume{34}(\bissue{2}),
\bfpage{683}--\blpage{713}
(\byear{2012})
\end{barticle}
\endbibitem

\bibitem[\protect\citeauthoryear{Goreinov et~al.}{2010}]{gostz-maxvol-2010}
\begin{bchapter}
\bauthor{\bsnm{Goreinov}, \binits{S.A.}},
\bauthor{\bsnm{Oseledets}, \binits{I.V.}},
\bauthor{\bsnm{Savostyanov}, \binits{D.V.}},
\bauthor{\bsnm{Tyrtyshnikov}, \binits{E.E.}},
\bauthor{\bsnm{Zamarashkin}, \binits{N.L.}}:
\bctitle{How to find a good submatrix}.
In: \beditor{\bsnm{Olshevsky}, \binits{V.}},
\beditor{\bsnm{Tyrtyshnikov}, \binits{E.}} (eds.)
\bbtitle{Matrix Methods: Theory, Algorithms, Applications},
pp. \bfpage{247}--\blpage{256}.
\bpublisher{World Scientific, Hackensack, NY},
\blocation{NY}
(\byear{2010})
\end{bchapter}
\endbibitem

\bibitem[\protect\citeauthoryear{Chen et~al.}{2021}]{chen2021solving}
\begin{barticle}
\bauthor{\bsnm{Chen}, \binits{Y.}},
\bauthor{\bsnm{Hosseini}, \binits{B.}},
\bauthor{\bsnm{Owhadi}, \binits{H.}},
\bauthor{\bsnm{Stuart}, \binits{A.M.}}:
\batitle{Solving and learning nonlinear pdes with gaussian processes}.
\bjtitle{Journal of Computational Physics}
\bvolume{447},
\bfpage{110668}
(\byear{2021})
\end{barticle}
\endbibitem

\bibitem[\protect\citeauthoryear{Meanti et~al.}{2022}]{meanti2022efficient}
\begin{bchapter}
\bauthor{\bsnm{Meanti}, \binits{G.}},
\bauthor{\bsnm{Carratino}, \binits{L.}},
\bauthor{\bsnm{De~Vito}, \binits{E.}},
\bauthor{\bsnm{Rosasco}, \binits{L.}}:
\bctitle{Efficient hyperparameter tuning for large scale kernel ridge
  regression}.
In: \bbtitle{International Conference on Artificial Intelligence and
  Statistics},
pp. \bfpage{6554}--\blpage{6572}
(\byear{2022})
\end{bchapter}
\endbibitem

\bibitem[\protect\citeauthoryear{Owhadi and Yoo}{2019}]{owhadi2019kernel}
\begin{barticle}
\bauthor{\bsnm{Owhadi}, \binits{H.}},
\bauthor{\bsnm{Yoo}, \binits{G.R.}}:
\batitle{Kernel flows: From learning kernels from data into the abyss}.
\bjtitle{Journal of Computational Physics}
\bvolume{389},
\bfpage{22}--\blpage{47}
(\byear{2019})
\end{barticle}
\endbibitem

\bibitem[\protect\citeauthoryear{Suykens}{2017}]{suykens2017deep}
\begin{barticle}
\bauthor{\bsnm{Suykens}, \binits{J.A.}}:
\batitle{Deep restricted kernel machines using conjugate feature duality}.
\bjtitle{Neural computation}
\bvolume{29}(\bissue{8}),
\bfpage{2123}--\blpage{2163}
(\byear{2017})
\end{barticle}
\endbibitem

\bibitem[\protect\citeauthoryear{Wenzel et~al.}{2023}]{wenzel2023data}
\begin{botherref}
\oauthor{\bsnm{Wenzel}, \binits{T.}},
\oauthor{\bsnm{Marchetti}, \binits{F.}},
\oauthor{\bsnm{Perracchione}, \binits{E.}}:
Data-driven kernel designs for optimized greedy schemes: A machine learning
  perspective.
arXiv preprint arXiv:2301.08047
(2023).
Accepted for publication in SISC.
\end{botherref}
\endbibitem

\bibitem[\protect\citeauthoryear{Narcowich et~al.}{2005}]{narcowich2005sobolev}
\begin{barticle}
\bauthor{\bsnm{Narcowich}, \binits{F.}},
\bauthor{\bsnm{Ward}, \binits{J.}},
\bauthor{\bsnm{Wendland}, \binits{H.}}:
\batitle{Sobolev bounds on functions with scattered zeros, with applications to
  radial basis function surface fitting}.
\bjtitle{Mathematics of Computation}
\bvolume{74}(\bissue{250}),
\bfpage{743}--\blpage{763}
(\byear{2005})
\end{barticle}
\endbibitem

\bibitem[\protect\citeauthoryear{Wendland and
  Rieger}{2005}]{wendland2005approximate}
\begin{barticle}
\bauthor{\bsnm{Wendland}, \binits{H.}},
\bauthor{\bsnm{Rieger}, \binits{C.}}:
\batitle{Approximate interpolation with applications to selecting smoothing
  parameters}.
\bjtitle{Numerische Mathematik}
\bvolume{101}(\bissue{4}),
\bfpage{729}--\blpage{748}
(\byear{2005})
\end{barticle}
\endbibitem

\bibitem[\protect\citeauthoryear{Wenzel et~al.}{2023}]{wenzel2023analysis}
\begin{barticle}
\bauthor{\bsnm{Wenzel}, \binits{T.}},
\bauthor{\bsnm{Santin}, \binits{G.}},
\bauthor{\bsnm{Haasdonk}, \binits{B.}}:
\batitle{Analysis of target data-dependent greedy kernel algorithms:
  Convergence rates for f-, f{\textperiodcentered} {P}-and f/{P}-greedy}.
\bjtitle{Constructive Approximation}
\bvolume{57}(\bissue{1}),
\bfpage{45}--\blpage{74}
(\byear{2023})
\end{barticle}
\endbibitem

\bibitem[\protect\citeauthoryear{Ma and Belkin}{2019}]{ma2019kernel}
\begin{barticle}
\bauthor{\bsnm{Ma}, \binits{S.}},
\bauthor{\bsnm{Belkin}, \binits{M.}}:
\batitle{Kernel machines that adapt to gpus for effective large batch
  training}.
\bjtitle{Proceedings of Machine Learning and Systems}
\bvolume{1},
\bfpage{360}--\blpage{373}
(\byear{2019})
\end{barticle}
\endbibitem

\bibitem[\protect\citeauthoryear{Goodfellow et~al.}{2016}]{goodfellow2016deep}
\begin{bbook}
\bauthor{\bsnm{Goodfellow}, \binits{I.}},
\bauthor{\bsnm{Bengio}, \binits{Y.}},
\bauthor{\bsnm{Courville}, \binits{A.}}:
\bbtitle{Deep Learning}.
\bpublisher{MIT Press},
\blocation{Cambridge}
(\byear{2016})
\end{bbook}
\endbibitem

\bibitem[\protect\citeauthoryear{Kingma and Ba}{2014}]{kingma2014adam}
\begin{botherref}
\oauthor{\bsnm{Kingma}, \binits{D.P.}},
\oauthor{\bsnm{Ba}, \binits{J.}}:
Adam: A method for stochastic optimization.
arXiv preprint arXiv:1412.6980
(2014)
\end{botherref}
\endbibitem

\bibitem[\protect\citeauthoryear{Telgarsky}{2015}]{telgarsky2015representation}
\begin{botherref}
\oauthor{\bsnm{Telgarsky}, \binits{M.}}:
Representation benefits of deep feedforward networks.
arXiv preprint arXiv:1509.08101
(2015)
\end{botherref}
\endbibitem

\bibitem[\protect\citeauthoryear{Jacot et~al.}{2018}]{jacot2018neural}
\begin{botherref}
\oauthor{\bsnm{Jacot}, \binits{A.}},
\oauthor{\bsnm{Gabriel}, \binits{F.}},
\oauthor{\bsnm{Hongler}, \binits{C.}}:
Neural tangent kernel: Convergence and generalization in neural networks.
{A}dvances in {N}eural {I}nformation {P}rocessing {S}ystems
\textbf{31}
(2018)
\end{botherref}
\endbibitem

\bibitem[\protect\citeauthoryear{Ehring and Haasdonk}{2023}]{ehring2023hermite}
\begin{botherref}
\oauthor{\bsnm{Ehring}, \binits{T.}},
\oauthor{\bsnm{Haasdonk}, \binits{B.}}:
Hermite kernel surrogates for the value function of high-dimensional nonlinear
  optimal control problems.
arXiv preprint arXiv:2305.06122
(2023)
\end{botherref}
\endbibitem

\bibitem[\protect\citeauthoryear{Dolgov et~al.}{2022}]{dolgov2022optimizing}
\begin{barticle}
\bauthor{\bsnm{Dolgov}, \binits{S.}},
\bauthor{\bsnm{Kalise}, \binits{D.}},
\bauthor{\bsnm{Saluzzi}, \binits{L.}}:
\batitle{Optimizing semilinear representations for state-dependent riccati
  equation-based feedback control}.
\bjtitle{IFAC-PapersOnLine}
\bvolume{55}(\bissue{30}),
\bfpage{510}--\blpage{515}
(\byear{2022})
\end{barticle}
\endbibitem

\bibitem[\protect\citeauthoryear{{Jones} and {Astolfi}}{2020}]{Astolfi2020}
\begin{bchapter}
\bauthor{\bsnm{{Jones}}, \binits{A.}},
\bauthor{\bsnm{{Astolfi}}, \binits{A.}}:
\bctitle{On the solution of optimal control problems using parameterized
  state-dependent {R}iccati equations}.
In: \bbtitle{2020 59th IEEE Conference on Decision and Control (CDC)},
pp. \bfpage{1098}--\blpage{1103}
(\byear{2020})
\end{bchapter}
\endbibitem

\end{thebibliography}

\end{document}